\documentclass{amsart}
\usepackage{amsfonts,amssymb,amscd,amsmath,enumerate,verbatim,calc}
\usepackage[all]{xy}

\newcommand{\CM}{Cohen-Macaulay}

\newcommand{\n}{\mathfrak{n} }
\newcommand{\m}{\mathfrak{m} }
\newcommand{\rr}{\mathfrak{r}}

\newcommand{\ZZ}{\mathbb{Z} }

\newcommand{\FF}{\mathbb{F}}

\newcommand{\Dc}{\mathcal{D} }

\newcommand{\X}{\mathbf{X} }

\newcommand{\rt}{\rightarrow}
\newcommand{\xar}{\longrightarrow}
\newcommand{\ov}{\overline}

\newcommand{\bu}{\mathbf{u}}

\newcommand{\bff}{\mathbf{f}}

\newcommand{\Pc}{\mathcal{P}}

\newcommand{\wt}{\widetilde }

\newcommand{\image}{\operatorname{image}}
\newcommand{\Tr}{\operatorname{Tr}}
\newcommand{\End}{\operatorname{End}}
\newcommand{\rad}{\operatorname{rad}}

\newcommand{\soc}{\operatorname{soc}}

\newcommand{\cx}{\operatorname{cx}}
\newcommand{\curv}{\operatorname{curv}}

\newcommand{\Syz}{\operatorname{Syz}}

\newcommand{\ann}{\operatorname{ann}}

\newcommand{\codim}{\operatorname{codim}}
\newcommand{\rank}{\operatorname{rank}}
\newcommand{\irr}{\operatorname{irr}}

\newcommand{\projdim}{\operatorname{projdim}}

\newcommand{\Hom}{\operatorname{Hom}}
\newcommand{\Ext}{\operatorname{Ext}}
\newcommand{\Tor}{\operatorname{Tor}}

\newcommand{\CMS}{\operatorname{\underline{CM}}}

\newcommand{\Gau}{\underline{\Gamma}}
\theoremstyle{plain}

\newtheorem{theorem}{Theorem}[section]
\newtheorem{corollary}[theorem]{Corollary}
\newtheorem{lemma}[theorem]{Lemma}
\newtheorem{proposition}[theorem]{Proposition}

\theoremstyle{definition}
\newtheorem{definition}[theorem]{Definition}
\newtheorem{s}[theorem]{}
\newtheorem{remark}[theorem]{Remark}

\theoremstyle{remark}

\begin{document}

\title{Symmetries and  connected components of the AR-quiver }
 \author{Tony J. Puthenpurakal}
\date{\today}
\address{Department of Mathematics, Indian Institute of Technology Bombay, Powai, Mumbai 400 076, India}
\email{tputhen@math.iitb.ac.in}
\subjclass{Primary 13 C14,  Secondary 13H10, 14B05}
\keywords{Artin-Reiten Quiver, Hensel rings, indecomposable modules, Ulrich modules, periodic modules, Non-periodic modules with bounded betti numbers}
\begin{abstract}
Let $(A,\m)$ be a commutative complete equicharacteristic  Gorenstein isolated singularity of dimension $d $ with 
$k = A/\m$  algebraically closed.   Let $\Gamma(A)$ be the AR (Auslander-Reiten) quiver of $A$.
Let  $\Pc$ be a property of maximal \CM \  $A$-modules.  We show that some naturally defined properties $\Pc$ define a union of connected components
of $\Gamma(A)$.   So in this case if there is a maximal \CM \ module satisfying $\Pc$  and if  $A$ is not of finite representation type  then there exists a family $\{ M_n \}_{n \geq 1}$ of maximal \CM \ indecomposable modules satisfying $\Pc$ with  multiplicity $e(M_n) > n$. Let $\Gau(A)$ be the stable quiver. We show that there are many symmetries in $\Gau(A)$. Furthermore $\Gau(A)$ is isomorphic to its reverse graph.  As an application we show that if $(A,\m)$ is a two dimensional Gorenstein isolated singularity with multiplicity $e(A) \geq 3$ then  for all $n \geq 1$ there exists an indecomposable self-dual maximal \CM \ $A$-module  of rank $n$. 
\end{abstract}

\maketitle

\section{introduction}
Let us recall that a commutative Noetherian local ring $(A,\m)$ is called an isolated singularity if $A_P$ is a regular local ring for all prime ideals $P \neq \m$. We note that with this definition if $A$ is Artinian then it is an isolated singularity. This is not a usual practice, nevertheless in this paper Artin rings will be considered as isolated singularities.  Also recall that if a local  Noetherian ring $(B,\n)$ is Henselian  then it satisfies Krull-Schmidt property, i.e., every finitely generated $B$-module is uniquely a direct sum of indecomposable $B$-modules.  Now assume that $B$ is \CM. Then we say $B$ is of finite  (\CM) representation type if $B$ has only finitely many indecomposable maximal \CM \ $B$-modules. In a remarkable paper Auslander proved that in this case $B$ is an isolated singularity, for instance see \cite[Theorem 4.22]{Y}.

To study
(not necessarily commutative) Artin algebra's
 Auslander and Reiten introduced the theory  of almost-split sequences.
These are now called AR-sequences. The AR-sequences are organized to form the AR-quiver. 
Later Auslander and Reiten extended the theory of AR-sequences to the case of commutative Henselian isolated singularities.

If $(A,\m)$ is a Henselian \CM \  isolated singularity then we denote its AR-quiver by $\Gamma(A)$.  
A good reference for this topic is \cite{Y}. 
The motivation for this paper comes from the following  crucial fact about AR quivers (under some conditions on $A$), see \cite[6.2]{Y}:

\text{If} \  $\mathcal{C}$  is a \emph{non-empty} connected component of $\Gamma(A)$ and if $A$ is not of finite representation type then there exist a family $\{ M_n \}_{n \geq 1}$ of maximal \CM \ indecomposable modules in $\mathcal{C}$ with  multiplicity $e(M_n) > n$.

Let $\Pc$ be a property of maximal \CM \ $A$-modules.   We show that some naturally defined properties $\Pc$ define a union of connected components of the AR quiver of $A$. 
Thus the above mentioned observation still holds. Therefore if there is a maximal \CM \ module satisfying $\Pc$ then there exists a family $\{ M_n \}_{n \geq 1}$ of maximal \CM \ indecomposable modules satisfying $\Pc$ with  multiplicity $e(M_n) > n$.

\s \label{hypothesis}  For the rest of the paper let us assume that $(A,\m)$  is a complete equicharacteristic Gorenstein isolated singularity of dimension $d$. 
 Assume $k = A/\m$ is algebraically closed.  Some of our results are applicable  more generally. However for simplicity  we will make this hypothesis throughout this paper. We will also assume that $A$ does not have finite representation type. This is automatic if $A$ is not a hypersurface, see \cite[8.15]{Y}. Furthermore if 
 $A$ is a hypersurface ring with $\dim A \geq 2$ and $e(A) \geq 3$ then also $A$ is not of finite representation type. 
 
  Now we describe our results. We first describe our results on 
  
  \emph{Connected Components of the AR-quiver:}

\textbf{I:}  \textit{Modules with periodic resolution:}\\
Let $M$ be a maximal \CM \  \emph{non-free}  $A$-module. Let $\Syz_n(M)$ be the $n^{th}$-syzygy module of $M$. We say $M$ has periodic resolution if there exists a non-negative integer $n$ and a positive integer $p$ with $\Syz_{n+p}(M) \cong \Syz_n(M)$.  
The smallest $p$ for which this holds is called the \emph{period} of the resolution.
We say $M$ has property  $\Pc$ if it has a periodic resolution. 

If $A$ is a hypersurface ring then any non-free maximal \CM \ $A$-module has periodic resolution with period $\leq 2$ and in fact $\Syz_3(M) \cong \Syz_1(M)$.
There exits maximal \CM \ modules with periodic resolutions if $A$ is a complete intersection of any codimension $c \geq 1$.  Again  it can be shown that in this case the period is $\leq 2$ and in fact $\Syz_3(M) \cong \Syz_1(M)$. 

For general Gorenstein local rings there is no convenient criterion to determine when $A$ has a module with perodic resolution (however see \cite[5.8]{C} for a criterion).
 It was conjectured by Eisenbud that if a module $M$ has a periodic resolution then the period is $\leq 2$, see \cite[p.\ 37]{E}. This was disproved by Gasharov and Peeva, see  \cite[Theorem 1.3]{GP}.

Our first result is
\begin{theorem}\label{first}
(with hypotheses as in  \ref{hypothesis}.) If $A$ is not a hypersurface ring then 
 $\Pc$ defines a union of connected components of  $\Gamma(A)$.
\end{theorem}

We now give more refined versions of Theorem \ref{first}:

\s \label{emb} Assume $A = Q/(f)$ where $(Q,\n)$ is 
a Gorenstein local ring and $f \in \n^2$ is a $Q$-regular element.  

Let  $M$ be a maximal \CM \ \emph{non-free} $A$-module. We say $M$ has property $\Pc_Q$ if   $\projdim_Q M$ finite. 
In this case   it is easy to prove that $M$ has a periodic-resolution over $A$  with period $\leq 2$. There is essentially a unique method to construct non-free maximal \CM \ modules over $A$  having finite projective dimension over $Q$. This is essentially due to Buchweitz \emph{et al,} see 
\cite[2.3]{BGS}.   Also see  the paper \cite[1.2]{HUB} by Herzog \emph{et al.}
Our next result is:

\begin{theorem}\label{em-comp}[ with hypotheses as in \ref{emb}.] If $Q$ is not regular then
$\Pc_Q$ defines a union of connected components of $\Gamma(A)$. 
\end{theorem}
Again  our results implies existence of indecomposable maximal \CM \ $A$-modules with arbitrarily   high multiplicity and satisfying property $\Pc_Q$.  However our method does not give a way to construct these modules.

\s  Eisenbud's conjecture (as stated above) is valid if $M$ has the so-called finite CI-dimension \cite[7.3]{AGP}.  We say $M$ has property $\Pc_O$  if  $M$ has finite CI-dimension over $A$ and has a periodic resolution over $A$. We say $M$ has property $\Pc_E$ if $M$ has periodic resolution over $A$ but  it has infinite CI dimension over $A$. Our next result is
\begin{theorem}\label{ci-comp}[with hypotheses as in \ref{hypothesis}.] Assume $A$ is not a hypersurface. 
\begin{enumerate}[\rm (1)]
\item
$\Pc_O$ defines a union of connected components  of $\Gamma(A)$. 
\item
$\Pc_E$ defines a union of connected components  of $\Gamma(A)$. 
\end{enumerate} 
\end{theorem} 
We note that  in \cite[3.1]{GP} a family $A_\alpha$ of an Artininian Gorenstein local ring  is constructed with each having a \emph{single}  module $M_\alpha$ having periodic resolution of period $> 2$ is given. 
As the period of $M_\alpha$ is greater than two it cannot have finite CI-dimension over $A_\alpha$.
Thus our result   
implies existence of indecomposable  modules with arbitrary length,  having a periodic resolution and having infinite CI dimension over $A_\alpha$. 

Note that till now  our results does not give any information regarding period's. In dimension two we can say something, see Theorem \ref{det-conn-comp}.

\s \label{ci} Now assume that $A$ is a complete intersection of codimension $c \geq 2$. There is a theory of support varieties for modules over $A$. Essentially  for every  finitely generated module $E$ over  $A$ an
 algebraic set $V(E)$ in the projective space $\mathbb{P}^{c-1}$ is attached, see \cite[6.2]{A}. Conversely it is known that if $V$ is an algebraic  set in $\mathbb{P}^{c-1}$ then there exists a finitely generated module $E$ with $V(E) = V$,
 see \cite[2.3]{B}. It is known that $V(\Syz_n(E)) = V(E)$ for any $n \geq 0$. Thus we can assume $E$ is maximal \CM. 
 If $E$ has periodic resolution over $A$ then $V(E)$ is  a finite set of points. The converse is also true. If further $E$ is indecomposable then $V(E)$ is a singleton set, see \cite[3.2]{B}. Let $a \in \mathbb{P}^{c-1}$. We say a maximal \CM \ $A$-module $M$ satisfies property  $\Pc_a$ if $V(M) = \{ a \}$. We prove:
 
 \begin{theorem}\label{ci-a}[ with hypotheses as in \ref{ci}.]
 Let $a \in \mathbb{P}^{c-1}$. Then $\Pc_a$ defines a union of connected components of $\Gamma(A)$.  Conversely if 
 $\mathcal{C}$ is a non-empty connected component of $\Gamma(A)$ containing a periodic module $M$ then for any $[N] \in \mathcal{C}$ we have $V(N) = V(M) (= \{ p \})$.  In particular $\Gamma(A)$ has at least  $|k|$ connected components.
 \end{theorem}

\textbf{II} \emph{Modules with bounded betti-numbers but not having a periodic resolution:}\\
For a long time it was believed that if a module $M$ has a bounded resolution (i.e., there exists $c$ with $\beta_i(M) \leq c$ for all $i \geq 0$)  then it is periodic.  If $A$ is a complete intersection then modules having bounded resolutions are periodic \cite[4.1]{E}. In \cite[3.2]{GP} there are examples of modules $M$ having a 
bounded resolution but $M$ is not periodic.

If $M$ is a maximal \CM \ $A$-module having a bounded resolution but $M$ is not periodic then we say that $M$ has property $\mathcal{B}_{NP}$.
We  prove 
\begin{theorem}\label{BNP}[with hypotheses as in \ref{hypothesis}.]
$\mathcal{B}_{NP}$ defines a union of connected components  of $\Gamma(A)$.
\end{theorem}
We note that if $M$ has bounded resolution but not periodic then there exists $c$ with $e(\Syz_n(M)) \leq c$ for all $n \geq 0$. Our result implies the existence of modules with bounded but not periodic resolution of arbitrary multiplicity.

\textbf{III:}  \emph{Ulrich modules:} \\
Let $M$ be a maximal \CM \ $A$-module. It is well-known that $e(M) \geq \mu(M)$ (here $\mu(M)$ denotes the cardinality  of a  minimal generating set of $M$) 
A maximal \CM \  module $M$ is said to be an Ulrich module if its multiplicity $e(M) = \mu(M)$. 
In this case we say $M$ has property $\mathcal{U}$. 

If $\dim A = 1$ then $A$ has a Ulrich module.  It is known that if $A$ is a strict complete intersection of any dimension $d$ then it has a Ulrich module, see \cite[2.5]{HUB}.  In particular if $A$ is a hypersurface ring then it has a Ulrich module. 
There are some broad class of  examples of Gorenstein normal domain (but not complete intersections) of dimension two that admit an Ulrich module see \cite[4.8]{BHU}. However there are no examples of Gorenstein local rings $R$  ( but not complete intersections) with $\dim R \geq 3$  such that $R$ admits an Ulrich module (note we are not even insisting that $R$ is reduced).

Even if $A$ is a hypersurface there is essentially a unique way to construct an Ulrich modules. We show

\begin{theorem}\label{u}[ with hypotheses as in \ref{hypothesis}.]
 Further assume that  either $A$ is a hypersurface ring of \emph{even} dimension $d \geq 2$ and multiplicity $e(A) \geq 3$   OR $A$ is Gorenstein of dimension two. 
Then $\mathcal{U}$ defines a union of connected components  of $\Gamma(A)$.
\end{theorem}

The reason we cannot say anything about Ulrich modules over hypersurface rings of \emph{odd} dimension 
is due to a peculiar nature of AR-sequences, see remark 8.2.  Also note that if $e(A) = 2$ then any non-free MCM $A$-module is an Ulrich module.

We now describe our result on:

\emph{Symmetries of AR-quiver:}\\
Let  $\Gamma_0(A)$ be the connected component of $\Gamma(A)$ containing the vertex $[A]$.  Set $\widetilde{\Gamma}(A) = \Gamma(A) \setminus \Gamma_0(A)$. 
Let $\Gau(A)$ denote the stable AR-quiver of $A$, i.e., we delete the vertex $[A]$ from $\Gamma(A)$ and all arrows connecting to $[A]$.  
Also set $\Gau_0(A)$ to be the stable part of $\Gamma_0(A)$.

Our starting point is the observation that for Klienian singularities $\Gau(A)$ is trivially isomorphic to its
reverse graph (see \cite[p.\  95]{Y}). Recall if $G$ is a directed graph then it's reverse graph $G^{rev}$ is a graph with the same vertices as $G$ and there is an arrow from vertex $u$ to $v$ in $G^{rev}$ if and only if there is an arrow from vertex $v$ to $u$ in $G$.  
In fact we construct
\begin{theorem}\label{rev} [with hypotheses as in \ref{hypothesis}.]
There exists isomorphisms  \\  $ D,  \lambda  \colon \Gau(A) \rt \Gau(A)^{rev}$ as graphs. If  $A$ is not a hypersurface ring then 
\begin{enumerate}[\rm (1)]
\item
$ D \neq \lambda$.
\item
There exists indecomposable maximal \CM \ modules  $M, N$ with \\  $\lambda(M) \neq M$ and $D(N) \neq N$. 
\end{enumerate}
\end{theorem}
The first isomorphism $D$ is just the dual functor.  The next map $\lambda$ arises in the theory of 
\emph{horizontal linkage} defined by Martsinkovsky and Strooker, see \cite[p.\ 592]{MS}.
We note that the assumption $A$ not a hypersurface is  essential for the later part of Theorem 
\ref{rev}, for in the case of Klienian singularities  it
is known that $\lambda(M) = M$ for each non-free indecomposable $M$, see \cite[Theorem 3]{MS}.

\s For $n \geq 0$ let $\Syz_n$ be the $n^{th}$ syzygy functor.
As $A$ is Gorenstein we can also define for integers $n \leq -1$ the 
$n^{th}$ cosyzygy  functor which we again denote with $\Syz_{n}$. 
By the definition of horizontal linkage we have $\Syz_{-1} \circ D = \lambda$. Thus $\Syz_{-1} = \lambda \circ D^{-1}$ and
$\Syz_1 = D \circ \lambda^{-1}$. So under the assumptions as in \ref{hypothesis} we get that $\Syz_n \colon  \Gau(A) \rt 
\Gau(A)$ is an isomorphism for all $n \in \ZZ$.
We prove:
\begin{theorem}\label{syzn}[ with hypotheses as in \ref{hypothesis}]
 Let $\mathcal{C}$  be a connected component of $\Gau(A)$.
For $[M] \in \mathcal{C}$, set $I(M) = \{ n \mid [\Syz_n(M)] \in \mathcal{C} \}$. Then
\begin{enumerate}[\rm(1)]
\item
$I(M)$ is an ideal in $\ZZ$ (possibly zero).
\item
$I(N) = I(M)$ for all $[N] \in \mathcal{C}$.
\end{enumerate}
\end{theorem} 

If $A$ is not of finite representation type then there is practically no information 
on connected components of $\Gau(A)$. The only case known is when  $A$ is a hypersurface there is information  on connected components of 
$\widetilde{\Gamma}(A)$, see \cite[Theorem I]{D}.
It is easy to show that $\Gau_0(A)$  has only finitely many components.
As an application of  Theorem \ref{syzn} we show:
\begin{corollary}\label{conn-0}
[with hypotheses as in \ref{hypothesis}.] Assume further that $A$ is not a hypersurface ring. Let $\mathcal{D}$ be a connected component of $\Gau_0(A)$. Then
\begin{enumerate}[\rm(1)]
\item
$\mathcal{D}$ has infinitely many vertices.
\item
The function $[M] \rt e(M)$ is unbounded on $\mathcal{D}$.
\end{enumerate}

\end{corollary}

\emph{Structure of the AR-quiver:}  \\
If $A$ is of finite representation type then the structure of the AR-quiver is known, see \cite{Y}.
For hypersurface rings which are not of finite representation type there is some information regarding connected components
of $A$ not containing the vertex $[A]$. For two dimensional Gorenstein rings we show:

\begin{theorem}\label{det-conn-comp}
[with hypothesis as in \ref{hypothesis}.]  Assume $\dim A = 2$ and $e(A) \geq 3$.
Let $\mathcal{C}$ be a non-empty component of $\Gamma(A)$. Then $\mathcal{C}$ is of the form
\[
M_1 \leftrightarrows M_2 \leftrightarrows M_3 \leftrightarrows  M_4 \leftrightarrows \cdots  \leftrightarrows M_n  \leftrightarrows \cdots
\]
where $e(M_n) = ne(M_1)$ for all $n \geq 1$. Furthermore 
\begin{enumerate}[\rm(1)]
\item
If $C = \Gamma_0(A)$ then $M_1 = A$. Furthermore $M_n^* \cong M_n$ for all $n \geq 1$. 
\item Assume now that $A$ is not a hypersurface ring. Then
\begin{enumerate}[\rm(a)]
\item
  If $M_j$ is periodic with period $c$ for some $j$ then  $M_n$ is periodic with period $c$ for all $n \geq 1$.  
  \item
  Let $\underline{\mathcal{C}}$ denote the stable part of $\mathcal{C}$. Let $[M_i] \in  \underline{\mathcal{C}}$. 
  If the Poincare series of $M_i$ is rational then the Poincare series of $M$ is rational for all $[M] \in  \underline{\mathcal{C}}$. Furthermore all of them  share a common denominator. 
\end{enumerate}
\end{enumerate}
\end{theorem}
In the Theorem above the Poincare series $P_M(z)$ of a module $M$ is \\ $\sum_{n \geq 0} \dim_k \Tor^A_n(M, k) z^n$.  
We also note that the structure of all components of $\Gamma(A) \setminus \Gamma_0(A)$ is already known, see \cite[4.16.2]{Ben}.

We have several interesting consequences of Theorem \ref{det-conn-comp}. A direct consequence of this theorem is that if 
$\dim A = 2$  and $e(A) \geq 3$ then for all $n \geq 1$ there exist's an indecomosable maximal \CM \ $A$-module $M_n$ of rank $n$ with $M_n^* 
\cong M_n$.
I do not know whether such a result holds for higher dimensional rings.

A simple consequence of Theorem's \ref{det-conn-comp} and \ref{syzn} is the following:
\begin{corollary}\label{d2-conn}[with hypotheses as in \ref{hypothesis}.] 
 Assume $A$ is not a hypersurface ring. Also assume $\dim A = 2$. Then  $\Syz_n(\Gau_0(A))$ are distinct for all $n \in \ZZ$, $n \neq 0$.
\end{corollary}
We now describe in brief the contents of this paper. In section two we discuss some preliminary results that we need.
In section three we discuss lifts of irreducible maps. In the next section we dicuss non-free indecomosable summands of maximal \CM \
approximation of the maximal ideal. In section five we give proof's of Theorem's 1.2, 1.4 and 1.6. In the next section we give
a proof of Theorem 1.8. In section seven we discuss our notion of quasi AR-sequences and in the next section we give a proof of 
Theorem 1.10.  In section nine we prove Theorem 1.11. In the next section we give  a few obstructions to existence of quasi-AR
sequences. In section 11 we describe the describe $\Gamma_0(A)$ when $A$ is a two dimensional with $e(A) \geq 3$.  In section twelve we give a proof of Theorem 1.15 and Corollary 1.16.
In the last section we discuss curvature and complexity of MCM modules and as an application give a proof of Theorem 1.9

\begin{remark}
 Srikanth Iyengar informed me about the excellent paper \cite{GZ} where the authors considered      AR-quiver of self-injective
 Artin algebra's. Note that commutative Artin Gorenstein rings is an extremely special case of self-injective Artin Algebra's.
 So our results in this case is sharper than that of \cite{GZ}. I do not beleive that the results of this paper when $A$ is
 commutative Artin Gorenstein ring
 will hold for the more general case of self-injective Artin algebra's. 
I also thank Dan Zacharia and Lucho Avramov for some useful discussions.
 \end{remark}

\section{Some Preliminaries}
In this paper all rings will be Noetherian local. All  modules considered are \textit{finitely} generated.
Let $(A,\m)$ be a local ring and let $k = A/\m$ be its residue field. If $M$ is an $A$-module then $\mu(M) = \dim_{k} M/\m M$ is the number of a minimal generating set of $M$. Also let $\ell(M)$ denote its length.  In this section we discuss a few preliminary results that we need.
 
\s \label{P} Let $M$ be an $A$-module. For $i \geq 0$  let $\beta_i(M) = \dim_k \Tor^A_i(M, k)$ be its $i^{th}$ \emph{betti}-number.  Let $P_M(z) = \sum_{n \geq 0}\beta_n(M)z^n$, the \emph{Poincare series} of $M$.
 Set 
 \[
 \cx(M) = \inf \{ d \mid \limsup \frac{\beta_n(M)}{n^{d-1}}  < \infty \}  \quad \text{and}  
 \]
 \[
 \curv{M} =  \limsup (\beta_n(M))^{\frac{1}{n}}.
 \]
 It is possible that $\cx(M) = \infty$, see \cite[4.2.2]{A-Inf}. However $\curv(M)$ is finite for any module $M$  \cite[4.1.5]{A-Inf}.  It can be shown that if $\cx(M) < \infty $ then $\curv(M) \leq 1$.  
 
 \s It can be shown that for any $A$-module $M$ we have
 \[
 \cx(M) \leq \cx(k) \quad \text{and} \quad \curv(M) \leq \curv(k). \ \ \ \  
 \]
 see \cite[4.2.4]{A-Inf}.
 
 \s \label{ext-ci} If $A$  is a complete intersection of co-dimension  $c$ then for any $A$-module $M$ we have $\cx(M) \leq c$. Furthermore for each $i = 0,\ldots,c$ there exists an $A$-module $M_i$ with $\cx(M_i) = i$. Also note that $\cx(k) = c$. \cite[8.1.1(2)]{A-Inf}. If  $A$ is a complete intersection and $M$ is a module with $\cx(M) = \cx(k)$ then we say $M$ is \emph{extremal}. 
 
 \s \label{ext-non-ci} If $A$ is \emph{not} a complete intersection then $\curv(k) > 1$. \cite[8.2.2]{A-Inf}. In this case we say a module $M$ is \emph{extremal} if $\curv(M) = \curv(k)$.

\s 
Let $M, N$ be $A$-modules and let $f \colon M \rt N$ be $A$-linear.  Let  $\FF_M$ be a minimal resolution of $M$ and let $\FF_N$ be a minimal resolution of $N$.  Then $f$ induces a lift $\widetilde{f} \colon \FF_M \rt \FF_N$.  This map $\widetilde{f}$  is unique up to  homotopy. The chain map $\widetilde{f}$ induces an $A$-linear map $f_n \colon \Syz_n(M) \rt \Syz_n(N)$ for all $n \geq 1$.

\s Denote by $\beta(M,N)$ the set of $A$-homomorphisms of $M$ to $N$ which pass through a  free module. That is, an $A$-linear map $f \colon M \rt N$ lies in $\beta(M, N)$ if and only if it factors as 
\[
\xymatrix{ 
M
\ar@{->}[d]_{u} 
\ar@{->}[dr]^{f}
 \\ 
F 
\ar@{->}[r]_{v} 
&N
}
\]
where $F \cong A^n$  for some $n \geq 1$.

\begin{remark}\label{lift-b}
If $f \colon M \rt N$ is $A$-linear and if $f_1, f_1^\prime \colon \Syz_1(M) \rt \Syz_1(N) $ are two lifts of $f$ then it is well known and easily verified that
$f_1 - f_1^\prime \in \beta(\Syz_1(M),\Syz_1(N))$.

\end{remark}

Recall that an $A$-module  $M$ is called \emph{stable} if $M$ has no free summands. 
We need the following:
\begin{proposition}\label{radE}
Let $(A,\m)$ be a Noetherian local ring and let $M, N$ be   $A$-modules. Set $\Lambda = \End_A(M)$ and  
$\rr = \rad \Lambda$.  Let $f \in \Lambda$. Also suppose there exists  $A$-linear maps $u \colon M \rt N$ and $v \colon N \rt M$. Set $g = v\circ u$.
\begin{enumerate}[ \rm(1)]
\item
If $f(M) \subseteq \m M$ then $f \in \rr$.
\item
If $M$ is stable and $f \in \beta(M,M)$ then $f \in \rr $.
\item
If $ 1 - g \in \rr$ then $u$ is a split mono and $v$ is a split epi.
\end{enumerate}
\end{proposition}
\begin{proof}
(1) This is well known. 

(2) Assume $f = \beta \circ \alpha$ where $\alpha \colon M \rt F$ and $\beta \colon F \rt M$ and $F  \cong A^n$. As $M$ is stable it follows that $\alpha(M) \subseteq \m F$ and so
 $f(M) \subseteq \m M$. The result follows from (1).
 
 (3) Let $1 - g = h$ where $h \in \rr$. Then $g = 1 - h$ is invertible in $\Lambda$.  So there exists $\tau \in \Lambda$ with $\tau \circ g = g \circ \tau = 1_M$. The result follows.
\end{proof}
\s Now assume $A$ is \CM. 
 Let $M, N$ be maximal \CM \ $A$-modules and let $f \colon M \rt N$ be $A$-linear.
Recall $f$ is said to be \emph{irreducible} if 
\begin{enumerate}
\item
$f$ is not a split epimorphism and not a split monomorphism.
\item
If $X$ is a maximal \CM \ $A$-module and if there  is a commutative diagram
\[
\xymatrix{ 
M
\ar@{->}[d]_{u} 
\ar@{->}[dr]^{f}
 \\ 
X 
\ar@{->}[r]_{v} 
&N
}
\]
 then either $u$ is a split monomorphism
or $v$ is a split epimorphism.
\end{enumerate}

\begin{remark}\label{syz-CM}
Let $(A,\m)$ be a \CM \ local ring. 
If $M$ is a maximal \CM \ $A$-module then it is easy to verify that $\Syz_n(M)$ is stable for each $n \geq 1$.  
\end{remark}

\s \label{dual} Suppose $M$ is a maximal \CM \ over a local Gorenstein ring $A$. Then $M^* = \Hom_A(M,A)$ is also a maximal \CM \ module. Furthermore $\Ext^i_A(M,A) = 0$ for $i > 0$.  We also have $(M^*)^* \cong M$.  Notice if
\[
0 \rt M_1 \xrightarrow{\alpha_1} M_2 \xrightarrow{\alpha_2} M_3 \rt 0,
\]
is a short exact sequence of maximal \CM \ $A$-modules then we have the following short exact sequence 
\[
0 \rt M_3^* \xrightarrow{\alpha_2^*} M_2^* \xrightarrow{\alpha_1^*} M_1^* \rt 0,
\]
of maximal \CM \ $A$-modules.

\begin{remark}\label{co-syzygy}
Let $(A,\m)$ be a Gorenstein local ring and let $M$ be a maximal \CM \ $A$-module. Then there exists exact sequences of the form
\[
0 \rt M \rt F \rt M_1 \rt 0,
\]
with $F$ free and $M_1$ maximal \CM.
\end{remark}

\begin{remark}\label{rem-dual} Due to \ref{dual} the following assertions hold:\\
Let $M_1, M_2$ and $N_1, N_2$ are maximal \CM \ $A$-modules and let $F, G$ be free $A$-modules. Suppose we have exact sequences:
\[
0 \rt M_1 \rt F \rt M_2 \rt 0 \quad \text{and} \quad 0 \rt N_1 \rt G \rt N_2 \rt 0. 
\]
If there exists an $A$-linear map $\psi_1 \colon M_1 \rt  N_1$ then:
\begin{enumerate}
\item
 there exists $A$-linear maps $\psi_2 \colon M_2 \rt N_2$ and $\phi \colon F \rt G$ such that the following diagram commutes:

\[
  \xymatrix
{
 0
 \ar@{->}[r]
  & M_1
    \ar@{->}[d]^{\psi_1}
\ar@{->}[r]
 & F
    \ar@{->}[d]^{\phi}
\ar@{->}[r]
& M_2
    \ar@{->}[d]^{\psi_2}
\ar@{->}[r]
 &0
 \\
 0
 \ar@{->}[r]
  & N_1
\ar@{->}[r]
 & G
\ar@{->}[r]
& N_2
\ar@{->}[r]
&0
 }
\]
\item
If $\psi^\prime \colon M_2 \rt N_2$ and $\phi^\prime \colon F \rt G$ are another pair of maps such that the above commutative diagram holds then $\psi_2 - \psi_2^\prime \in \beta(M_2, N_2)$. 
\end{enumerate}
\end{remark}
\begin{definition}\label{pre-lift}
(with hypotheses as in \ref{rem-dual}.) We call $\psi_2$ to be a \emph{pre-lift} of $\psi_1$.
\end{definition}

The following  is an easy consequence of \ref{dual}.
\begin{proposition}\label{dual-irr}
Let $(A,\m)$ be a local Gorenstein ring and let $M,N$ be maximal \CM \ $A$-modules. Let $f \colon M \rt N$ be $A$-linear and let $f^* \colon N^* \rt M^*$ be the induced map.
We have:
\begin{enumerate}[\rm (1)]
\item
$f$ is a split mono if and only if $f^*$ is a split epi.
\item
$f$ is a split epi if and only if $f^*$ is a split mono.
\item
$f$ is irreducible if and only if $f^*$ is irreducible. 
\end{enumerate} \qed
\end{proposition}

\s Let $(A,\m)$ be a Gorenstein local ring.  
Let us  recall the definition of Cohen-Macaulay approximation from \cite{AB}. A  \CM \ approximation of an
  $A$-module $M$ is an exact sequence
\[
0 \xar Y \xar X \xar M \xar 0,
\]
where $X$ is a maximal \CM \ $A$-module and $Y$ has finite projective dimension.
Such a sequence is not unique but $X$ is known to unique up to a free summand and so is well defined in the stable category   $\CMS(A)$. We denote by $X(M)$ the maximal \CM \ approximation of $M$

\s \label{codim-n}  If $M$ is \CM \ then maximal \CM \ approximation of $M$ are very easy to construct. We recall this construction from \cite{AB}.  Let $n = \codim M = \dim A - \dim M$.  Let $M^\vee = \Ext^n_A(M,A)$. It is well-known that $M^\vee$ is \CM \ module of codim $n$ and $M^{\vee \vee} \cong M$.
Let $\FF$ be any free resolution of $M^\vee$ with each $\FF_i$ a finitely
generated free module. Note $\FF$ need not be minimal free resolution of $M$. Set  $S_n(\FF) = \image( \FF_n \xrightarrow{\partial_n} \FF_{n-1})$. Then  note $S_n(\FF)$ is maximal \CM \ $A$-module. It can be easily proved that $X(M) \cong S_n(\FF)^*$ in $\CMS(A)$.

The following result is well-known and easy to prove.
\begin{proposition}\label{syz-mcm-approx}
Let $(A,\m)$ be a Gorenstein local ring and let $M$ be an $A$-module. Then
$$   X(\Syz_1(M)) \cong \Syz_1(X(M)) \quad \text{in} \ \CMS(A). $$
\end{proposition}

\section{Lifts of irreducible maps}
In this section $(A,\m)$ is a Gorenstein local ring, not necessarily an isolated singularity.  Also $A$ need not be Henselian. 

 The following is the main result of this section.

\begin{theorem}\label{lift}
(with hypotheses as above.) Let $M, N$ be stable maximal \CM \ $A$-modules and let $f \colon M \rt N$ be $A$-linear. Let $f_1 \colon \Syz_1(M) \rt \Syz_1(N)$ be any lift of $f$.
If $f$ is irreducible then $f_1$ is also an irreducible map.
\end{theorem}

We need a few  preliminaries to prove Theorem \ref{lift}.

We first prove:
\begin{lemma}\label{add-beta}
Let $(A,\m)$ be a local Gorenstein ring and let $M, N$ be stable maximal \CM \ $A$-modules. Let $f \colon M \rt N$ be $A$-linear and let $\delta \in \beta(M,N)$.
If $f$ is irreducibe then so is $f + \delta$. 
\end{lemma}
\begin{proof}
Assume $\delta = v \circ u$ where $u \colon M \rt F$ and $v \colon F \rt N$ and $F  \cong A^n$. As $M$ is stable it follows that $u(M) \subseteq \m F$ and so
 $\delta(M) \subseteq \m N$. 
 
 \textit{Claim-1:} $f + \delta $ is not a split mono. \\
 Suppose it is so. Then there exists $\sigma \colon N \rt M$ with $\sigma \circ (f + \delta) = 1_M$. So $\sigma \circ f  + \sigma \circ \delta = 1_M$.
 As $\delta(M) \subseteq \m N$ we get $\sigma \circ \delta(M) \subseteq \m M$. It  follows that $\sigma \circ \delta \in \rad \End_A(M)$.
   By \ref{radE}(3) it  follows that $f$ is a split mono. This is  a contradiction as $f$ is irreducible.
 
 \textit{Claim-2:} $f + \delta$ is not a split epi. \\
 Suppose it is so. Then there exists $\sigma \colon N \rt M$ with $(f + \delta)\circ \sigma =  1_N$. So $f\circ \sigma + \delta \circ \sigma = 1_N$.  Notice 
 $$\delta \circ \sigma(N)  \subset \delta(M) \subseteq \m N.$$
 It follows that $\delta \circ \sigma \in \rad \End_A(N)$.  By \ref{radE}(3)  it follows that $f$ is a split epimorphism. This is a contradiction as $f$ is irreducible. 
 
 \textit{Claim 3:}  Suppose $X$ is maximal \CM \ and  we have a commutative diagram 
 \[
\xymatrix{ 
M
\ar@{->}[d]_{g} 
\ar@{->}[dr]^{f + \delta}
 \\ 
X 
\ar@{->}[r]_{h} 
&N
}
\]
 then either $g$ is a split monomorphism
or $h$ is a split epimorphism. \\
\textit{Proof of Claim 3:} Notice we have a commutative diagram
 \[
\xymatrix{ 
M
\ar@{->}[d]_{(g, -u)} 
\ar@{->}[dr]^{f }
 \\ 
X\oplus F
\ar@{->}[r]_{h+v} 
&N
}
\]
As $f$ is irreducible either $(g,-u)$ is a split mono or $h + v$ is a split epi. We  assert: \\
\textit{Subclaim-1:} If $(g,-u)$ is a split mono then $g$ is a split mono. \\
\textit{Subclaim-2:} If $h + v$ is a split epi then $h$ is a split epi.\\
 Notice that Subclaim 1 and 2  will finish the proof of Claim 3. Also Claims 1,2,3  implies the assertion of the Lemma.
 We now give:\\
 \textit{Proof of Subclaim-1:} As $(g,-u)$ is a split mono there exits  $\sigma \colon X \oplus F \rt M$ such that
 $\sigma \circ( g, -u) = 1_M$. Write $\sigma  = \sigma_1 + \sigma_2$ where  $\sigma_1 \colon X \rt M$ and $\sigma_2 \colon F \rt M$.  Thus we have 
 $\sigma_1\circ g - \sigma_2 \circ u = 1_M$. As $M$ is stable $u(M) \subseteq \m F$. So $\sigma_2\circ u(M) \subseteq \m M$. Thus $\sigma_2\circ u \in \rad \End_A(M)$.
 It follows  from \ref{radE}(3)
 that $g$ is a split mono.

 We now give:\\
 \textit{Proof of  Subclaim-2:}
 As $h + v$ is a split epi there exists $\sigma \colon N \rt X \oplus F$ with  $(h+v) \circ \sigma = 1_N$. Write $\sigma = (\sigma_1, \sigma_2)$ where $\sigma_! \colon N \rt X$ and $\sigma_2 \colon N \rt F$. It follows that $h\circ \sigma_1 + v \circ \sigma_2 = 1_N$.
 
 As $N$ is stable $\sigma_2(N) \subseteq \m F$. So $v \circ \sigma_2(N) \subseteq \m N$. Thus $v \circ \sigma_2 \in \rad \End_A(N)$. It follows  from \ref{radE}(3) that $h$ is a split epi.
 \end{proof}

We also need
\begin{lemma}\label{lift-mono}
Let $(A,\m)$ be a Gorenstein local ring and let $M$ be a stable maximal \CM \ $A$-module. Let $F = A^{\mu(M)}$ and let $\epsilon \colon F \rt M$ be a minimal map. Set $M_1 = \ker \epsilon \cong \Syz_1(M)$.  Let $X$ be another maximal \CM \  $A$-module (not necessarily stable) and let $\eta \colon G \rt X$ be a projective cover (not necessarily minimal). Let $X_1 = \ker \eta$. Let $\alpha \colon M \rt X$ be $A$-linear and let $\alpha_1 \colon M_1 \rt X_1$ be any lift of $\alpha$.
If $\alpha$ is a split mono then $\alpha_1$ is a split mono.
\end{lemma}
\begin{proof}
We note that $M_1$ is also stable.  Let $\phi \colon X \rt M$ be such that $\phi \circ \alpha = 1_M$. Let $\phi_1 \colon X_1 \rt M_1$ be a lift of $\phi$. Then note that $\phi_1 \circ \alpha_1$ is a lift of $1_M$. Thus $\phi_1\circ \alpha_1 - 1 \in \beta(M_1,M_1)$. The result now follows from \ref{radE}.
\end{proof}

The following is a dual version of Lemma \ref{lift-mono} and can be proved similarly.
\begin{lemma}\label{lift-epi}
Let $(A,\m)$ be a Gorenstein local ring and let $N$ be a stable maximal \CM \ $A$-module. Let $F = A^{\mu(N)}$ and let $\epsilon \colon F \rt N$ be a minimal map. Set $N_1 = \ker \epsilon \cong \Syz_1(N)$.  Let $X$ be another maximal \CM \  $A$-module (not necessarily stable) and let $\eta \colon G \rt X$ be a projective cover (not necessarily minimal). Let $X_1 = \ker \eta$. Let $\beta \colon X \rt N$ be $A$-linear and let $\beta_1 \colon X_1 \rt N_1$ be any lift of $\alpha$.
If $\beta$ is a split epi then $\beta_1$ is a split epi. \qed
\end{lemma}

We now give:
\begin{proof}[Proof of Theorem \ref{lift}.]
Set $M_1 = \Syz_1(M)$ and $N_1 = \Syz_1(N)$.\\
\textit{Claim-1:} $f_1$ is not a split mono. \\
Suppose if possible $f_1$ is a split mono. Then there exists $\sigma_1 \colon N_1 \rt M_1$ with $\sigma_1 \circ f_1  = 1_{M_1}$. Let $\sigma \colon N \rt M$ be a pre-lift of $\sigma$ (see \ref{pre-lift} for this notion). Then $\sigma \circ f$ is a pre-lift of $1_{M_1}$. Notice $1_M$ is a pre-lift of of $1_{M_1}$. Then by \ref{rem-dual} we get that $1_{M} - \sigma \circ f \in \beta(M,M)$. As $M$ is stable we get by \ref{radE} that $f$ is a split mono. This is a contradiction as $f$ is irreducible.

\textit{Claim-2:} $f_1$ is not a split epi. \\
Suppose if possible $f_1$ is a split epi. Then there exists $\sigma_1 \colon M_1 \rt N_1$ with $  f_1 \circ \sigma_1 = 1_{N_1}$. Let $\sigma \colon N \rt M$ be a pre-lift of $\sigma$.  Then $ f \circ \sigma $ is a pre-lift of $1_{N_1}$. Notice $1_N$ is a pre-lift of of $1_{N_1}$. Then by \ref{rem-dual} we get that $1_{N} - f \circ \sigma \in \beta(N,N)$. As $N$ is stable we get by \ref{radE} that $f$ is a split epi. This is a contradiction as $f$ is irreducible.

\textit{Claim-3:}
If $X_1$ is a maximal \CM \ $A$-module and if there  is a commutative diagram
\[
\xymatrix{ 
M_1
\ar@{->}[d]_{u_1} 
\ar@{->}[dr]^{f_1}
 \\ 
X_1 
\ar@{->}[r]_{v_1} 
&N_1
}
\]
 then either $u_1$ is a split monomorphism
or $v_1$ is a split epimorphism.\\
\textit{Proof of Claim-3:} By \ref{co-syzygy} there exists an exact sequence
\[
0 \rt X_1 \rt L \rt X \rt 0,
\]
with $L_1$ free and $X$ maximal \CM. 

Let $u \colon M \rt X$ be a pre-lift of $u_1$ and let $v \colon X \rt N$ be a pre-lift of $v_1$. Then notice $v \circ u$ is a pre-lift of $f_1 = v_1 \circ u_1$. As $f$ by definition is a pre-lift of $f_1$ we get that $v\circ u = f + \delta $ for some $\delta \in \beta(M, N)$.  

By \ref{add-beta} we get that $f + \delta $ is irreducible. So $u$ is a split mono or $v$ is a split epi. By Lemma's \ref{lift-mono} and \ref{lift-epi} we get that $u_1$ is a split mono or $v_1$ is a split epi.

By Claims 1, 2 and 3 the result follows. 
\end{proof}

\section{Indecomposable non-free summands of \\ maximal \CM \ approximation of the maximal ideal}
In this section $(A,\m)$ is a Henselian Gorenstein local ring. Let $X(\m)$ be a maximal \CM \ approximation of the maximal ideal. In this section we are concerned with non-free indecomposable summands of $X(\m)$. Our results are:

\begin{theorem}\label{mcm-extremal}
Let $(A,\m)$ be a Henselian Gorenstein local ring of dimension $d$ and  let $X(\m)$ be a maximal approximation of $\m$. Let $M$ be an indecomposable non-free summand of $X(\m)$. Then $M$ is extremal, i.e.,
\begin{enumerate}[\rm (1)]
\item
If $A$ is a complete intersection of codimension $c$ then $\cx(M) = \cx(k) = c$.
\item
If $A$ is not a complete intersection then $\curv(M) = \curv(k) >1$.
\end{enumerate}
\end{theorem}

We also prove:
\begin{theorem}\label{mcm-Ulrich}
Let $(A,\m)$ be a Henselian Gorenstein local ring of dimension $d \geq 1$ and infinite residue field $k$. Let $e(A) \geq 3$. Assume either $\dim A = 2$ or $A$ is a hypersurface ring (with no restriction on dimension) with multiplicity $e(A) \geq 3$. Let $M$ be an indecomposable Ulrich $A$-module. Then neither $M$ or $\Syz_1(M)$ is a summand of $X(\m)$. 
\end{theorem}

\s \emph{Motivation:} Our motivation to prove the above results is the following: Assume $A$ is a Gorenstein Henselian isolated singularity. If $M$ is a maximal \CM \  \emph{non-free} indecomposable module then there exists an irreducible morphism from $M \rt A$  only if $M$  is a summand of $X(\m)$, see \cite[4.2.1]{Y}. 
In our Theorems  we have to show that the vertex $[A]$ does not belong to certain components of $\Gamma(A)$.

We first give:
\begin{proof}[Proof of Theorem \ref{mcm-extremal}]
By Proposition \ref{syz-mcm-approx} we get that $X(\m) \oplus F  =  \Syz_1(X(k))\oplus G$ for some free modules $F,G$. 
Thus it suffices to prove that if $M$ is a direct summand of $X(k)$ then it is extemal.  By \ref{codim-n} it suffices to prove that if $M$ is a summand of $\Syz_d(k)^*$  then $M$ is extremal.  We prove by induction on $d$ that if $M$ is a summand of $\Syz_n(k)^*$ for some $n \geq d$ then $M$ is extremal.

We first consider the case $d = 0$. We note that as $k$ is indecomposable $\Syz_n(k)$ is indecomposable for all $n \geq 0$. So $\Syz_n(k)^*$ is indecomposable. Therefore $M = \Syz_n(k)^*$. Notice $\Syz_n(\Syz_n(k)^*) = k$. It follows that $M$ is extremal. 

We now assume that $d \geq 1$ and the result has been proved for Gorenstein Henselian rings of dimension $d -1$.
Let $x \in \m \setminus \m^2$ be a non-zero divisor on $A$.  Set $B = A/(x)$ and for any $A$-module $N$ set $\ov{N} = N/xN$.
We note that 
\[
\ov{\Syz^A_d(k)} \cong  \Syz^B_d(k) \oplus \Syz^B_{d-1}(k).
\]
It follows that
\[
\ov{\Syz^A_d(k)^*} \cong \Syz^B_d(k)^* \oplus \Syz^B_{d-1}(k)^*.
\]
If $M$ is a summand of $\Syz^A_d(k)^*$ then $\ov{M}$ is a summand of  $\ov{\Syz^A_d{k}^*}$. Let $E$ be an irreducible summand of $\ov{M}$. Then by Krull-Schmidt
it is an irreducible summand of either $\Syz^B_d(k)^*$ or of $\Syz^B_{d-1}(B)^*$. By induction hypothesis  we get that $E$ is extremal. It follows that $M$ is extremal.

We prove by induction on $n \geq d$ that if $M$ is an irreducible summand of $\Syz^A_n(k)^*$ then $M$ is extremal. We just proved the result for $n = d$.
We now assume that $n \geq d + 1$ and the result has been proved for $n-1$.  There is an exact sequence
\[
0 \rt \Syz_n^A(k) \rt F \rt \Syz^A_{n-1}(k) \rt 0, \quad \text{with} \ F \ \text{free.}
\]
As $n \geq d + 1$ we have that $\Syz_{j}^A(k)$ is maximal \CM \ for all $j \geq n-1$. So there is an exact sequence
\[
0 \rt \Syz_{n-1}^A(k)^* \rt F^* \rt \Syz_n^A(k)^* \rt 0.
\]
Therefore $\Syz_1^A(M)$ is a summand of $\Syz_{n-1}(k)^*$. By induction hypothesis $\Syz_1^A(M)$ is extremal. So $M$ is extremal.
\end{proof}

We now give:
\begin{proof}[Proof of Theorem \ref{mcm-Ulrich}]
\emph{Case 1:} We first consider the case when $A$ is a hypersurface of dimension $d \geq 1$ and multiplicity $e(A) \geq 3$. 

Let $M$ be an indecomposable Ulrich A-module. Suppose if possible  $M$ is a summand 
of $X(\m)$. By Proposition \ref{syz-mcm-approx} we get that $X(\m) \oplus F  =  \Syz_1(X(k))\oplus G$ for some free modules $F,G$.  It follows that $\Syz_{-1}(M)$ is a summand of $X(k)$. As $M$ has no free summands we get $\Syz_1(M) = \Syz_{-1}(M)$. By \ref{codim-n} we get  $\Syz_1(M)$ is a summand of $\Syz_d(k)^*$. It follows that $\Syz_1(M)^*$  is a summand of $\Syz_d(k)$. But 
$$\Syz_1(M)^* \cong \Syz_{-1}(M^*) \cong \Syz_{1}(M^*).$$
Notice if $M$ is Ulrich then $M^*$ is also Ulrich.  Similarly if we set $N = \Syz_{1}(M)$ is a summand of $\X(\m)$ then $M^*$ is a summand of $\Syz_d(k)$. 

By the arguments in the previous paragraph   it suffices to prove that if $E$ is an Ulrich $A$-module then neither $E$ nor $\Syz_{1}(E)$  is a summand of $\Syz_d(k)$. This we prove by induction on $d$.

We first consider the case $d = 1$. Then as $e(A) \geq 3$ we have that $\m = \Syz_1(k)$ is indecomposable \cite[Theorem A]{T}.

 If $E$ is a summand of $\m$ then $\m = E$.
Let $x$ be $E$-superficial. Then as $E = \m$ is Ulrich we get that $\m \m = x \m$. So $\m^2 = x \m$. So $A$ has minimal multiplicity. It follows that $e(A) = 2$. This is a contradiction.

  If $\Syz_1(E)$ is a summand of $\m$ then $\m = \Syz_1(E)$. Using 
  \cite[Theorem 2]{P} we get the $h$-polynomial of $\m$ 
  \[
  h_\m(z) = 2(1 + z + z^2 + \cdots + z^{e-2}) \quad \text{where} \ e = e(A) \geq 3.
  \]
  It follows that the $h$-polynomial of $A$ is
  \[
  h_A(z) = 1 + z(h_\m(z) - 1)  =  2z^{e-1} + \text{lower terms in $z$}. 
  \]
  This is a contradiction as $A$ is a hypersurface. 

Now assume that $d \geq 2$ and the result has been proved for hypersurface rings of dimension $d-1$.
Let $x \in \m \setminus \m^2$ be sufficiently general. Then $x$ is $A$-regular and $A \oplus E \oplus \Syz^A_1(E)$-superficial. Set $B = A/(x)$ and if $V$ is an $A$-module set $\ov{V} = V/xV$. Then notice 
\[
\ov{\Syz_d^A(k)} = \Syz_{d}^B(k) \oplus \Syz_{d-1}^B(k).
\]
Note $\ov{E}$ is an Ulrich $B$-module. Let $\ov{E} = U_1 \oplus U_2 \oplus \cdots \oplus U_s $ where $U_i$ are indecomposable $B$-modules. Then each $U_i$ is an Ulrich $B$-module. We also have 
$$ \ov{\Syz_1^A(E)} \cong \Syz^B_1(\ov{E}) \cong \Syz^B_1(U_1) \oplus\cdots \oplus \Syz^B_1(U_s). $$
 By \cite[8.17]{Y}, $\Syz^B_1(U_i)$ is an indecomposable $B$-module for $i = 1,\ldots, s$. 
 
 If $E$ is a summand of $\Syz_d^A(k)$ then $\ov{E}$ is a summand of 
 $\ov{\Syz_d^A(k)}$. So $U_1$ is a summand of $\Syz_{d}^B(k)$ or  $\Syz_{d-1}^B(k)$.
 By our induction hypothesis $U_1$ is not a summand of $\Syz_{d-1}^B(k)$. It follows that $U_1$ is a summand of $\Syz_{d}^B(k)$. Therefore $\Syz^B_1(U_1)$ is a summand of $\Syz_{d+1}^B(k) \cong \Syz_{d-1}^B(k)$, a contradiction. A similar argument will show that $\Syz_1^A(E)$ is not a summand of $\Syz_d^A(k)$.

 \emph{Case 2:} We now consider the case when $A$ is a Gorenstein local ring  of dimension $2$ but not a hypersurface.\\
 Notice  $e(A) \geq 4$.  Let $M$ be an indecomposable non-free, maximal \CM \ $A$-module. 
 If $M$ is a summand of  $X(\m)$ then $\Syz_{-1}^A(M)$ is a summand of $X(k)$. Notice $\Syz_{-1}^A(M)$ is an
 indecomposable non-free $A$-module. Therefore $\Syz_{-1}^A(M)$ is a summand of $\Syz_2^A(k)^*$. 
 By \cite[Theorem B]{T}, $\Syz_2^A(k)$ is an indecomposable maximal \CM \ $A$-module. 
 It follows that $\Syz_{-1}^A(M) \cong \Syz_2^A(k)^*$. Notice 
 \begin{equation}\label{U-Gor}
 \Syz_1(M^*) \cong \Syz_{-1}^A(M)^*  \cong \Syz_2^A(k).
 \end{equation}
Let $x_1,x_2$ be a $A \oplus M^* \oplus \Syz_1(M^*)$-superficial sequence. Set $C = A/(x_1, x_2)$ and
if $E$ is an $A$-module, set $\ov{E} = E/(x_1, x_2)$. We note that
\begin{equation}\label{U-Gor-2}
\ov{\Syz_2^A(k)} = \Syz_2^C(k) \oplus \Syz_1^C(k)^2 \oplus \Syz_0^C(k).
\end{equation}

 Set $N = \ov{M^*}$. Then $\ov{\Syz_1^A(M^*)} \cong \Syz_1^C(N)$. We now consider two cases.\\
 \emph{Case 1:} $M$ is Ulrich.\\
 Then notice $M^*$ is also Ulrich. Then $N = k^a$ where $a = \mu(N)$.  Let $\n$ be the maximal ideal of $C$. Note it is indecomposable as a $C$-module. By \ref{U-Gor-2} and \ref{U-Gor} we get
 \[
 \n^a = \Syz_2^C(k) \oplus \Syz_1^C(k)^2 \oplus k.
 \]
By Krull-Schmidt we get that $k = \n$. So $\n^2 = 0$. It follows that $e(A) = e(C) = 2$, a contradiction.

\emph{Case 2:} $M = \Syz_1^A(D)$  where $D$ is Ulrich.\\
Then $D = \Syz_{-1}^A(M)$. It follows that $\Syz_2^A(k)^*$ is Ulrich. Therefore $\Syz_2^A(k)$
is also Ulrich. Therefore $\ov{\Syz_2^A(k)} \cong k^l$ for some $l$. So by $\ref{U-Gor-2}$ it
follows that $\n$ the maximal ideal  of $C$ is isomorphic to $k$. Therefore $\n^2 = 0$ and so $e(A) = e(C) = 2$ 
a contradiction.
\end{proof}

\section{Proofs of Theorem's \ref{first}, \ref{em-comp} and \ref{ci-comp}}
We first give
\begin{proof}[Proof of Theorem \ref{first}]
Let $M$ be an indecomposable periodic maximal \\ Cohen-Macaulay $A$-module. Let $N, L$ be indecomposable  non-free maximal \CM  \ $A$-modules and assume there exists irreducible maps $u \colon N \rt M$ and $v \colon M \rt L$.  Let
\[
0 \rt \tau(M) \rt E_M \rt M \rt 0, \quad 0 \rt M \rt V_M \rt \tau^{-1}(M) \rt 0
\]
be the AR-sequences starting and ending at $M$. 

We first note that  $\Syz_i(u) \colon  \Syz_i(N) \rt \Syz_i(M)$  and $\Syz_i(v) \colon \Syz_i(M) \rt \Syz_i(L)$ are irreducible, see \ref{lift}. Let $p$ be the period of $M$. Note $\Syz_p(M) \cong M$.

We have irreducible maps $\Syz_{ip} \colon \Syz_{ip}(N) \rt M$ for all $i \geq 1$.
As $\Syz_{ip}(N)$ is indecomposable maximal \CM \ $A$-module we get that $\Syz_{ip}(N)$ are factors of $E_M$ for all $i \geq 1$. By Krull-Schmidt theorem we get that $N$ is periodic. A dual argument gives that $L$ is periodic.

If there exists an irreducible map
from $M \rt A$ then $M$ is factor of $X(\m)$ the maximal \CM \ approximation of $\m$.  
So by \ref{mcm-extremal} we get that $M$ is extremal. As $A$ is not a hypersurface this is a contradiction.

Notice $\tau(M) = \Syz_{-d + 2}(M)$ as $A$ is Gorenstein. If there is an irreducible map from $A \rt M$ then $A$ is a factor of $E_M$ and so  there exist's an irreducible map from $\tau(M) \rt A$. By previous argument we get that $\tau(M)$ is extremal. So $M$ is extremal. As $A$ is not a hypersurface this is a contradiction,

Thus $\Pc$ defines a union of connected components of $\Gamma(A)$. 
\end{proof}
\begin{remark}
 In Proposition 5.2 of \cite{AR} it is shown that if $A$ is a self-injective Artin algebra  and $M, N$ are non-projective $A$-modules
 such that there is an irreducible map $f \colon M \rt N$  then $M$ is a periodic module if and only if $N$ is. To the
 best of my knolwedge this proof does not generalize to higher dimensions.
\end{remark}

We now give 
\begin{proof}[Proof of Theorem \ref{em-comp}]
Let $M$ be an indecomposable maximal \CM \ $A$-module such that $\projdim_Q M$ is finite. Then $M$ is a periodic $A$-module with period $\leq 2$.  As $A$ is not a hypersurface ring then by 
proof of previous Theorem there is no irreducible map from $A \rt M$ or $M \rt A$. 

 Let $N, L$ be indecomposable  non-free maximal \CM  \ $A$-modules and assume there exists irreducible maps $u \colon N \rt M$ and $v \colon M \rt L$.  Let
\[
0 \rt \tau(M) \rt E_M \rt M \rt 0, \quad 0 \rt M \rt V_M \rt \tau^{-1}(M) \rt 0
\]
be the AR-sequences starting and ending at $M$.  Notice $\tau(M) = \Syz_{-d + 2}^A(M)$ and $\tau^{-1}(M) = \Syz_{d -2}^A(M)$.
It follows that $\projdim_Q \tau(M)$ and $\projdim_Q \tau^{-1}(M)$ is finite.
Therefore $\projdim_Q E_M$ and $\projdim_Q V_M$ is finite. 

 By \cite[5.5, 5.6]{Y}, we get that  $N, L$ are direct summands of $E_M$ and $V_M$ respectively. It follows that $\projdim_Q N$ and $\projdim_Q L$ are finite. Thus $\Pc_Q$ determines a union of connected components of $\Gamma(A)$.
\end{proof}

\s Let us recall that a quasi-deformation $A \rt B <-- Q$ of $A$ is a flat local map $A \rt B$ and a deformation  $Q \xrightarrow{\eta} B$ (i.e., $\ker \eta$ is generated by a $Q$-regular sequence. 
We say
 CI-dimension of an $A$-module $M$ is finite if there is a quasi-deformation $A \rt B <--Q$ with $\projdim_Q M\otimes_A B$ is finite. 
 We say CI-dimension of $M$ is infinite if it is not finite.
 
 We now give:
 \begin{proof}[Proof of Theorem \ref{ci-comp}]
 Let $M$ be an indecomposable maximal \CM \ $A$-module which is periodic and has a finite CI-dimension over $A$. 
 Let  $A \rt B <--Q$ be a quasi-deformation of $A$ with $\projdim_Q M\otimes_A B$ finite.
 As $A$ is not a hypersurface  by  proof of  Theorem 1.2 there is no irreducible map from $A \rt M$ or $M \rt A$. 

 Let $N, L$ be indecomposable  non-free maximal \CM  \ $A$-modules and assume there exists irreducible maps $u \colon N \rt M$ and $v \colon M \rt L$.  Let
\[
0 \rt \tau(M) \rt E_M \rt M \rt 0, \quad 0 \rt M \rt V_M \rt \tau^{-1}(M) \rt 0
\]
be the AR-sequences starting and ending at $M$.  Notice $\tau(M) = \Syz_{-d + 2}^A(M)$ and $\tau^{-1}(M) = \Syz_{d -2}^A(M)$. 

Notice  $\tau(M)\otimes_A B =  \Syz_{-d+2}^B(M\otimes_A B)$. Thus CI-dimension of $\tau(M)$ is finite. As the quasi-deformations involved are the same we get that $E_M$ also has  
finite CI-dimension over $A$. A similar argument yields that $V_M$ has finite CI dimension over $A$. As $N$ and $L$ are summands of $E_M$ and $V_M$ respectively we get that CI dimension of $N$ and $L$ are finite.
It follows that  $\Pc_O$ defines a union of connected components of $\Gamma(A)$.

Let $\mathcal{C}$ be the union of connected components of $\Gamma(A)$ consisting of periodic indecomposable  maximal \CM \ $A$-modules and let $\mathcal{C}_O$ be the union of connected components of $\Gamma(A)$ consisting of periodic indecomposable maximal \CM \ $A$-modules having finite CI-dimension over $A$. 
Then as $\mathcal{C}_O \subseteq \mathcal{C}$ we get that $\mathcal{C} \setminus \mathcal{C}_O$ is a union of connected components of $\Gamma(A)$.
Notice $\mathcal{C} \setminus \mathcal{C}_O$
consists of precisely  those periodic maximal \CM \ $A$-modules  which has infinite CI-dimension over $A$.
 \end{proof}

\section{Proof of Theorem \ref{ci-a}}
We need to recall some preliminaries regarding support varieties. This is relatively simple in our case since $A$ is complete with algebraically closed residue field.

\s Let $A = Q/(\bu)$ where $(Q,\n)$ is regular local and $\bu = u_1,\ldots, u_c \in \n^2$ is a regular sequence.
We need the notion of cohomological operators over a complete intersection ring; see \cite{Gull} and
\cite{E}.
The \emph{Eisenbud operators}, \cite{E}  are constructed as follows: \\
Let $\mathbb{F} \colon \cdots \rightarrow F_{i+2} \xrightarrow{\partial} F_{i+1} \xrightarrow{\partial} F_i \rightarrow \cdots$ be a complex of free
$A$-modules.

\emph{Step 1:} Choose a sequence of free $Q$-modules $\wt{F}_i$ and maps $\wt{\partial}$ between them:
\[
\wt{\mathbb{F}} \colon \cdots \rightarrow \wt{F}_{i+2} \xrightarrow{\wt{\partial}} \wt{F}_{i+1} \xrightarrow{\wt{\partial}} \wt{F}_i \rightarrow \cdots
\]
so that $\mathbb{F} = A\otimes\wt{\mathbb{F}}$

\emph{Step 2:} Since $\wt{\partial}^2 \equiv 0 \ \text{modulo} \ (\mathbf{u})$, we may write  $\wt{\partial}^2  = \sum_{j= 1}^{c} u_j\wt{t}_j$ where
$\wt{t_j} \colon \wt{F}_i \rightarrow \wt{F}_{i-2}$ are linear maps for every $i$.

 \emph{Step 3:}
Define, for $j = 1,\ldots,c$ the map $t_j = t_j(Q, \mathbf{f},\mathbb{F}) \colon \mathbb{F} \rightarrow \mathbb{F}(-2)$ by $t_j = A\otimes\wt{t}_j$.

\s
The operators $t_1,\ldots,t_c$ are called Eisenbud's operator's (associated to $\mathbf{u}$) .  It can be shown that
\begin{enumerate}
\item
$t_i$ are uniquely determined up to homotopy.
\item
$t_i, t_j$ commute up to homotopy.
\end{enumerate}
\s Let $R = A[t_1,\ldots,t_c]$ be a polynomial ring over $A$ with variables $t_1,\ldots,t_c$ of degree $2$. Let $M, N$ be  finitely generated $A$-modules. By considering a free resolution $\mathbb{F}$ of $M$ we get well defined maps
\[
t_j \colon \Ext^{n}_{A}(M,N) \rightarrow \Ext^{n+2}_{R}(M,N) \quad \ \text{for} \ 1 \leq j \leq c  \ \text{and all} \  n,
\]
which turn $\Ext_A^*(M,N) = \bigoplus_{i \geq 0} \Ext^i_A(M,N)$ into a module over $R$. Furthermore these structure depend  on $\bu$, are natural in both module arguments and commute with the connecting maps induced by short exact sequences.

\s  Gulliksen, \cite[3.1]{Gull},  proved that 
$\Ext_A^*(M,N) $ is a finitely generated $R$-module. We note that $\Ext^*(M,k)$ is a finitely generated graded module over $T = k[t_1,\ldots, t_c]$. Define
$V^*(M) = Var(\ann_T(\Ext^*(M,k))$ in the projective space $\mathbb{P}^{c-1}$. We call 
$V^*(M)$ the support variety of a module $M$.

We need the following 
\begin{lemma}\label{existence-indec-irred}
Let $(Q,\n)$ be a complete regular local ring  with  algebraically closed residue field $k$.  Let  $\bff = f_1,\ldots, f_c \in \n^2$ be a regular sequence. Assume $c \geq 2$. Set $A = Q/(\bff)$   and let $d = \dim A$.  Let $W$ be an irreducible non-empty sub-variety of $\mathbb{P}^{c-1}$. Then there exists an indecomposable non-free maximal \CM \ $A$-module $M$ with $V^*(M) = W$.  
\end{lemma}
\begin{proof}
By \cite[2.3]{B},   there exists an $A$-module $E$ with $V^*(E) = W$. Then $N = \Syz_{d+1}(E)$ is a maximal \CM \  $A$-module and $V^*(N) = V^*(E) =  W$. Notice $N$ has no free summands.   If $N$ is indecomposable then we are done. Otherwise $N = N_1 \oplus N_2$ where $N_1$ and $N_2$ are maximal \CM \ $A$-modules with no free-summands.
Let $T = k[t_1,\ldots,t_c]$ and let $\Ext^*_A(N,k)$, $\Ext^*_A(N_1,k)$ and $\Ext^*_A(N_2,k)$ be given $T$-module structure as above. 
As $\Ext^*(N,k) = \Ext^*_A(N_1,k) \oplus \Ext^*_A(N_2,k)$ we get that
\[
\ann_T \Ext^*(N,k) = \ann_T \Ext^*(N_1,k) \cap \ann_T \Ext^*(N_2,k).
\]
It follows that 
\[
W = W_1 \cap W_2 \quad \text{where} \ W_i = V^*(N_i) \ \text{for} \ i = 1,2.
\]
As $W$ is irreducible we get that $W = W_1$ or $W = W_2$. Iterating this procedure we get our result.
\end{proof}

The following result yields Theorem \ref{ci-a} as an easy corollary.
\begin{theorem}\label{conn-irr}
Let $Q = k[[x_1,\ldots,x_n]]$ be the formal power series over an algebraically closed field $k$. Let  $\bu = u_1,\ldots, u_c \in \n^2$ be a regular sequence. Assume $c \geq 2$. Set $A = Q/(\bu)$   and let $d = \dim A$.  Let $W$ be a proper sub-variety  of $\mathbb{P}^{c-1}$. Let 
\[
\mathcal{C}_W = \{ M \mid M \  \text{is indecomposable MCM \  $A$-module with $V^*(M) = W$ }\}.
\]
Then $\mathcal{C}_W$ defines an union of connected components of $\Gamma(A)$. If $W$ is irreducible then $\mathcal{C}_W$ is non-empty.
\end{theorem}
\begin{proof}
Let $M \in  \mathcal{C}_W$. Notice $M$ is not free. As $W$ is a proper subset of $\mathbb{P}^{c-1}$ we get  that $\dim W \leq c-2$. So $\cx M  = \dim W + 1 \leq c-1$. In particular $M$ is not extremal. So there is no irreducible map $M \rt A$. As $\tau(M) = \Syz_{-d+2}(M)$ is also not extremal there is no irreducible map $\tau(M) \rt A$. So there is no irreducible map $A \rt M$.

Let $N, L$ be  indecomposable non-free maximal \CM \ $A$-modules and suppose there exists an irreducible map $u \colon N \rt M$ and an irreducible map $v \colon M \rt L$.\\
\textit{Claim:} $V^*(N) = V^*(M) = W$ and $V^*(L) = V^*(M) = W$. \\
Suppose there exists $a \in V^*(N) \setminus V^*(M)$. Let 
 $D$ be a maximal \CM \ $A$-module with $V^*(D) = \{ a \}$. As $\{ a \} \cap W = \emptyset$ we get that $\Ext^A_i(D, M) = 0$ for all $i \gg 0$. As $V^*(\tau(M)) = V^*(M) = W$ we also get $\Ext^A_i(D,\tau(M)) = 0$ for $i \gg 0$. Thus $\Ext^A_i(D, E_M) = 0$ for $i \gg 0$. As $N$ is a summand of $E_M$ we get that  $\Ext^A_i(D, N) = 0$ for all $i \gg 0$. This implies $a \notin V^*(N)$, a contradiction. Thus $V^*(N) \subseteq V^*(M)$.
 We now notice that there exists an irreducible map $M \rt \tau^{-1}(N)$. By the previous argument we get that  $V^*(M) \subseteq V^*(\tau^{-1}(N)) = V^*(N)$. Thus $V^*(N) = V^*(M) = W$.
 
 As there exist's an irreducible map $v \colon M \rt N$. Then there exists an irreducible map from $v^\prime \colon  N \rt \tau^{-1}(M)$. By  the previous argument we get
 $V^*(N) = V^*(\tau^{-1}(M))$. As $V^*(\tau^{-1}(M)) = V^*(M) = W$ we get $V^*(N) = W$. Thus we have proved our Claim. 
 
 By our claim and as there are no irreducible maps from $M \rt A$ and $A \rt M$ we get that $\mathcal{C}_W$ is a union of connected components of $\Gamma(A)$. 
 If $W$ is irreducible then by \ref{existence-indec-irred} we get that $\mathcal{C}_W$ is non-empty.
\end{proof}
We now give
\begin{proof}[Proof of Theorem \ref{ci-a}]
By \ref{conn-irr} the result follows.
\end{proof}

\section{quasi AR-sequences}
\s \label{setup-quasi}\textit{Setup:} In this section $(A,\m)$ is a Henselian Gorenstein local ring with algebraically closed residue field $k$. We also assume $A$ is an isolated singularity.  For the notion of AR -sequences see \cite[Chapter 2]{Y}. In this section we introduce the notion of  quasi AR-sequences.

\begin{definition}
Let $M$ be an indecomposable non-free maximal \CM \ $A$-module. By a \emph{quasi-AR sequence} \textbf{ending} at $M$ we mean an exact sequence $s \colon 0 \rt K \rt E \xrightarrow{\phi} M $ such that  
\begin{enumerate}
\item
$E$ is a stable maximal \CM \ $A$-module.
\item
$\phi$ is irreducible.
\item
If $L$ is a stable  maximal \CM \ $A$-module and if there is an $A$-linear map $\sigma \colon L \rt M$  which is not a split epi then there exist's a map $\xi \colon L \rt E$ such that
$\phi \circ \xi - \sigma \in \beta(L,M)$.
\end{enumerate}
\end{definition}
\begin{remark}
Unlike AR-sequences the module $K$ \textit{need not} be maximal \CM. Also the map $\phi$ need not be surjective.
\end{remark}

A consequence of the definition of quasi AR-sequence is the following:
\begin{proposition}\label{qar-conseq}[with hypothesis as in \ref{setup-quasi}.] 
Suppose $\sigma$ is irreducible. Then $\xi$ is a split monomorphism.
\end{proposition}
\begin{proof}
By \ref{add-beta},  $\phi \circ \xi$ is irreducible. Now $\phi$ is irreducible. So in particular it is not a split epi. It follows that $\xi$ is a split mono.
\end{proof}

We need the following analogue to Corollary 2.12 from \cite{Y}.
\begin{lemma}\label{dir-qar}
Let  $(A,\m)$ be a Henselian Gorenstein local ring with algebraically closed residue field $k$. Let $M, L$ be indecomposable non-free maximal \\ \CM \ $A$-modules and let $s \colon 0 \rt K \rt E \xrightarrow{\phi} M $ be a quasi AR-sequence ending at $M$. Then the following two conditions are equivalent:
\begin{enumerate}[\rm (i)]
\item
There is an irreducible morphism from $L$ to $M$.
\item
$L$ is isomorphic to a direct summand of $E$.
\end{enumerate}
\end{lemma}
\begin{proof}
$(i) \implies (ii)$. This follows from \ref{qar-conseq}.

$(ii) \implies (i)$. Assume the decomposition of $E$ is given by $E = L \oplus Q$. Denote $\phi = (f,g)$ along this decomposition. We claim that $f$ is irreducible.

Clearly $f$ is not a split epi as $\phi$ is not a split epi.
If $f$ is a split mono then as $L$ and $M$ are indecomposable we get that $f$ is an isomorphism and so a split epi, a contradiction.

The rest of the proof is similar to the proof of $(ii) \implies (i)$ of Corollary 2.12 from \cite{Y}.
\end{proof}

We give two constructions of quasi AR-sequences. The first one comes from AR sequences.
\begin{proposition}\label{ar-qar}
Let  $(A,\m)$ be a Henselian Gorenstein local ring with algebraically closed residue field $k$. Let $M$ be an indecomposable non-free maximal \CM \ $A$-modules and let $l \colon 0 \rt N \rt E_M \xrightarrow{p} M  \rt 0$ be an  AR-sequence ending at $M$.
Suppose $ E_M = E\oplus F$ where $F$ is a free $A$-module and $E$ has no free summands. Denote $p = (\phi, \psi)$ along this decomposition. Assume $E \neq 0$. Let $K = \ker \phi$.
Then $s \colon 0 \rt K \rt E \xrightarrow{\phi} M $ is  a quasi AR-sequence ending at $M$.
\end{proposition}
\begin{proof}
By Corollary 2.12 from \cite{Y} we get that $\phi$ is an irreducible map. Now let $L$ be a stable maximal \CM \ $A$-module and let $f \colon L \rt M$ be $A$-linear which is not a split epi. Then as $l$ is an AR sequence ending at $M$ there exist's an $A$-linear map $g \colon L \rt E_M$ with $f = p \circ g$.  Suppose
\[
g = \left(\begin{matrix} \sigma \\ \delta \end{matrix} \right) \quad \text{where} \ \sigma \colon L \rt E \ \text{and} \ \delta \colon L \rt F.
\]
So we get $f = \phi \circ \sigma + \psi \circ \delta$. Notice 
$\psi \circ \delta \in \beta(L,M)$. It follows that $s$ is a quasi AR sequence ending at $M$. 
\end{proof}

\s Let $l \colon 0 \rt N \rt E_M \xrightarrow{p} M  \rt 0$ be an  AR-sequence ending at $M$ then it is not true that a lift of $p$;  $q \colon \Syz_1(E) \rt \Syz_1(M)$ is surjective and defines  a AR sequence ending at $\Syz_1(M)$. The great advantage of quasi AR sequences is that it behaves well under lifting (and also pre-lifting).

\begin{theorem}\label{lift-qar}
Let  $(A,\m)$ be a Henselian Gorenstein local ring with algebraically closed residue field $k$. Let $M$ be an indecomposable non-free maximal \CM \ $A$-modules and let $s \colon 0 \rt K \rt E \xrightarrow{\phi} M $ be a quasi AR-sequence ending at $M$. 
Let $\psi$ be any lift of $\phi$.  Set $K^\prime = \ker \psi$.
Then $s^\prime \colon 0 \rt K^\prime \rt \Syz_1(E) \xrightarrow{\psi} \Syz_1(M)$ is a quasi AR sequence ending at $\Syz_1(M)$.  Similarly if $\theta$ is any pre-lift of $\phi$. Then $\widetilde{s} \colon  0 \rt \widetilde{K}  \rt \Syz_{-1}(E) \xrightarrow{\theta} \Syz_{-1}(M)$ is a quasi AR sequence ending at $\Syz_{-1}(M)$.
\end{theorem}
\begin{proof}
As $E$ is a stable maximal \CM \ $A$-module we get that $\Syz_1(E)$ is also a stable maximal \CM \ $A$-module. 
By Theorem \ref{lift} we get that $\psi$ is an irreducible map. 

Let $L$ be a stable maximal \CM \ $A$-module and let $f \colon L \rt \Syz_1(M)$ be an $A$-linear map which is not a split epi. Let $g \colon \Syz_{-1}(L) \rt M$ be any pre-lift of $f$. Then by \ref{lift-mono} we get that $g$ is not a split epi. It follows that there exists $\xi \colon \Syz_{-1}(L) \rt M$ such that
$\phi \circ \xi - g  = \delta$ where $\delta \in \beta(\Syz_{-1}(L), M)$.
Let $\xi^\prime  \colon L \rt \Syz_1(M)$ be a lift of $\xi$. Then notice by construction $\psi \circ \xi^\prime  - g$ is a lift of $\delta$. It follows that
$\psi \circ \xi^\prime  - g \in \beta(L, \Syz_1(M))$. Thus $s^\prime$ is a quasi AR-sequence ending at $\Syz_1(M)$.

The assertion regarding $\widetilde{s}$ can be proved similarly. 
\end{proof}

\s Till now we have not used the fact that  $k$, the residue field of $A$ is algebraically closed. We will now use this fact.
Let us recall the following definition from \cite[Chapter 5]{Y}. \\ 
Let $M, N$ be maximal \CM  \ $A$-modules. Set $(M,N) = \Hom_A( M, N)$. 

Decompose $M = \bigoplus_{i = 1}^{m} M_i$ and $N = \bigoplus_{ j = 1}^{n} N_j$ where $M_i, N_j$ are indecomposable $A$-modules for all $i, j$
For $g \in (M, N)$ decompose $g = (g_{ij})$ where $g_{ij} \colon M_i \rt N_j$. 

\begin{definition}
We say $g \in (M, N)_* $ if no $g_{ij}$ is an isomorphism.
\end{definition}

\s \label{irr} We define the following descending chain $\{ (M, N)_n \}_{n \geq 1} $ of $A$-submodules of $(M,N)$ as follows:

$(M,N)_n$ consists of those $f \in (M,N)$ such that there is a sequence $X_0,\ldots, X_n$ of maximal \CM  \ $A$-modules with $X_0 = M$ and $X_n = N$  and $f_i \in (X_{i-1}, X_i)_*$  such that $f = f_n \circ f_{n-1} \circ \cdots \circ f_1$.

It is easy to see that $(M,N)_n$ are $A$-submodules of $(M, N)$ and that
\[
(M, N) \supseteq (M, N)_1 \supseteq  \cdots \supseteq (M, N)_n \supseteq (M, N)_{n+1} \supseteq \cdots.
\]
It is not difficult to see that $(M,N)_1/ (M,N)_2$ is a $k = A/\m$ vector space. It is finite dimensional since it is finitely generated as an $A$-module. Set
\[
\irr(M, N) = \dim_k \frac{(M,N)_1}{(M,N)_2}.
\] 

Let us restate the following basic result from \cite[5.5]{Y}.
\begin{lemma}\label{basic-Y}[with hypothesis as in \ref{setup-quasi}.] 
 Let $M, N$ be indecomposable maximal \CM \ $A$-modules. Assume there is an AR-sequence ending at $M$
$$ 0 \rt \tau(M)  \rt E_M \rt M \rt 0. $$
Let $n$ be the number of copies  of $N$ in direct summands of $E_M$ (note that $n = 0$ is possible). Then the following equality holds:
$$ \irr(N, M) =  n.$$
\end{lemma} 

We note that  the assumption $k$ is algebraically closed is used in the proof of Lemma \ref{basic-Y}.  The following is a basic result in our theory of quasi AR-sequences.
\begin{theorem}\label{basic-qar}[with hypothesis as in \ref{setup-quasi}.] 
 Let $M, N$ be indecomposable \emph{non-free}  maximal \CM \   $A$-modules. Let $0 \rt K \rt E \xrightarrow{\phi} M$ be a quasi  AR sequence ending at $M$.
Let $n$ be the number of copies of $N$ in direct summands of $E$ (note that $n = 0$ is possible). Then the following equality holds:
$$ \irr(N, M) = n.$$
\end{theorem}
\begin{proof}
Set $S(N,E) = (N,E)/(N,E)_1$. Then by proof of Lemma 5.5. in \cite{Y} it follows that $S(N,E) \cong k^n$. 
Define
\begin{align*}
\theta \colon S(N,E) &\rt  \frac{(N,M)_1}{(N,M)_2}, \\
[f ] &\rt  [\phi \circ f].
\end{align*}
By  \ref{dir-qar}  we get that $\theta$ is a well-defined $k$-linear map. 

We first show that $\theta$ is surjective.  Let $\sigma \colon N \rt M$ be an irreducible map. Denote by $[\sigma]$ it's class in $(N,M)_1/(N,M)_2$. By our definition of quasi AR sequence there exist's  $\xi \colon N \rt E$ such that  $\phi \circ \xi - \sigma  =  g \in \beta(N, M)$.  By \ref{qar-conseq} we get that $\xi$ is a split monomorphism. As $N, M$ are both non-free indecomposable $A$-modules we get that $g \in (N,M)_2$. Therefore we get 
$ \theta([\xi])  = [\sigma]$. 

Next we show that $\theta$ is injective
. Let $[h] \in S(N,E)$ be non-zero. Thus $h \colon N \rt E$ is a split mono. Then by \ref{dir-qar} we get that $ \phi \circ h \colon N \rt M$ is
 an irreducible map. Thus $\theta([h]) = [\phi \circ h] \neq 0$. Therefore $\theta$ is injective. 
 \end{proof}
 
A consequence of the previous two results is the following:
\begin{corollary}\label{qar-ar}
[with hypothesis as in \ref{setup-quasi}.]   Let $M$ be indecomposable \emph{non-free}  maximal \CM \   $A$-module. Suppose the following is an AR-sequence ending at $M$:
$$ t: \ 0 \rt \tau(M) \rt E_M \xrightarrow{p}  M \rt 0. $$
Further assume $E_M$ has no free summnads. 
Let $$s: 0 \rt K \rt E \xrightarrow{\phi} M$$ be a quasi  AR sequence ending at $M$. Then $E \cong E_M$ and $\phi$ is surjective.  Furthermore $s$ is also an AR-sequence ending at $M$ (and so $K \cong \tau(M)$). 
\end{corollary} 
\begin{proof}
By \ref{basic-Y} and \ref{basic-qar} it follows that $E \cong E_M$. As $p \colon E_M \rt M$ is an indecomosable map, by defining property of quasi AR-sequences there exists a  map $\xi \colon E_M \rt E$  such that $\phi \circ \xi - p = \delta \in \beta(E_M, M)$. As $p$ is indecomposable, by \ref{qar-conseq} we get $\xi$ is a split mono. As $E \cong E_M$,   by Krull-Schmidt we get $\xi$ is an isomorphism.
It follows that
\[
\phi - p \circ \xi^{-1} = \delta \circ \xi^{-1} := \eta \in \beta(E,M).
\] 
Set $\psi = p \circ \xi^{-1} \colon E \rt M$. Notice $\psi$ is surjective. 
As $E, M$ has no free summands we get that $\eta(E) \subseteq \m M$. It follows that the maps $\ov{\phi}, \ov{\psi} \colon E/\m E \rt M/\m M$ are equal. As $\ov{\psi}$ is surjective, it follows that $\ov{\phi}$ is surjective. So by Nakayama's Lemma we get that $\phi$ is surjective. 

We now use that  $t$ is an AR-sequence. As $\phi$ is irreducible, it is not a split epi. Therefore there exists $\theta \colon E \rt E_M$ such that $p \circ \theta = \phi$. As $\phi$ is irreducible we get that $\theta$ is a split mono. Since $E \cong E_M$,   by Krull-Schmidt we get $\theta$ is an isomorphism. Note there exists $f \colon  K \rt \tau(M)$ which makes the following diagram commute:
\[
  \xymatrix
{
s \colon
 0
 \ar@{->}[r]
  & K
    \ar@{->}[d]^{f}
\ar@{->}[r]
 & E
    \ar@{->}[d]^{\theta}
\ar@{->}[r]^{\phi}
& M
    \ar@{->}[d]^{1_M}
\ar@{->}[r]
 &0
 \\
 t \colon
 0
 \ar@{->}[r]
  & \tau(M)
\ar@{->}[r]
 & E_M
\ar@{->}[r]^{p}
& M
\ar@{->}[r]
&0
 }
\]
By Snake Lemma we get that $f$ is an isomorphism. So $K \cong \tau(M)$.
In the terminology of \cite[2.3]{Y} we get $s \sim t$. So $s$ is an AR-sequence ending at $M$.
\end{proof} 
The following consequence of Corollary \ref{qar-ar} is significant:
\begin{lemma}\label{min-gen}[with hypothesis as in \ref{setup-quasi}.] 
 Let $M$ be an indecomposable maximal \CM \  non-free $A$-module. Let $t \colon 0 \rt \tau(M) \rt E_M \xrightarrow{p} M \rt 0$ be an AR-sequence ending at $M$. If there is no irreducible maps $A \rt M$ and $A \rt \Syz_1(M)$ then we have
$$ \mu(E_M) = \mu(M) + \mu(\tau(M)). $$
\end{lemma}
\begin{proof}
Set $N = \tau(M), M_1 = \Syz_1(M)$ and $E_1 = \Syz_1(E)$. As there is no irreducible maps from $A \rt M$ we get that $E_M$ is stable. In particular $t$ is a quasi AR-sequence, see \ref{ar-qar}. Let $\phi \colon E_1 \rt M_1$ be any lift of $p$. Then 
we have a quasi AR-sequence 
\[
s \colon  0 \rt K_1 \rt E_1 \xrightarrow{\phi} M_1.
\]
As there are no irreducible maps from $A \rt M_1$ we get that $E_{M_1}$ is stable. Therefore by Corollary \ref{qar-ar} we get that $\phi$ is surjective and $s$ is an AR-sequence ending at $M_1$. Furthermore $E_1 \cong E_{M_1}$ and $K_1 \cong \tau(M_1)$.

Let $F \rt E_M$ and $G \rt M$ be projective covers. Then we have an exact sequence

\[
  \xymatrix
{
 0
 \ar@{->}[r]
  & E_1
    \ar@{->}[d]^{\phi}
\ar@{->}[r]
 & F
    \ar@{->}[d]^{\theta}
\ar@{->}[r]
& E_M
    \ar@{->}[d]^{p}
\ar@{->}[r]
 &0
 \\
 0
 \ar@{->}[r]
  & M_1
\ar@{->}[r]
 & G
\ar@{->}[r]
& M
\ar@{->}[r]
&0
 }
\]
Set $N = \tau(M)$.
As $\phi, p$ are surjective we get $\theta$ is surjective. As $G$ is free it is in fact a split epi. Set $H = \ker \theta$. Then $H$ is free and  $\mu(H) = \mu(E_M) - \mu(M)$.  

As $\phi$ is surjective we get by Snake Lemma that the induced map $H \rt N$ is surjective.
It follows that  $\mu(N) \leq  \mu(E_M) - \mu(M)$.  As there is an exact sequence
 $0 \rt N \rt E_M \rt M \rt 0$ it follows (after tensoring with $A/\m$) that
  $\mu(N) \geq  \mu(E_M) - \mu(M)$. The result follows.
\end{proof}

\section{Proof of Theorem \ref{u}}
In this section we give a proof of Theorem \ref{u}.

\s\label{even-dim} [with hypotheses as in \ref{hypothesis}.]  Recall if $M$ is an indecomposable non-free $A$-module then
it's AR-translate is $\tau(M) = \Syz_{-d + 2}(M)$. So if $d = 2$ or if $A$ is a hypersurface of even dimension then 
$\tau(M) = M$.

We now give:
\begin{proof}[Proof of Theorem \ref{u}.]
By \ref{even-dim} we get that $\tau(M) = M$. Let $M$ be an indecomosable Ulrich $A$-module. Then by \ref{mcm-Ulrich} 
there is no irreducible map from $ M \rt A$ and $\Syz_1(M) \rt A$. By our assumptions on the ring there is no irreducible 
map from $A \rt M$ and $A \rt \Syz_1(M)$. Let $s \colon 0 \rt M \rt E_M \rt M \rt 0$ be the AR-sequence ending at $M$. Then
by \ref{min-gen} we get $\mu(E_M) = 2\mu(M)$. Also note that $e(E_M) = 2 e(M)$. So $e(E_M) = \mu(E_M)$. Therefore $E_M$ is Ulrich 
$A$-module.

Let $N$ be a non-free indecomposable maximal \CM \ $A$-module. If there is an irreducible morphism $N \rt M$ then $N$ is a summand
of $E_M$. As $E_M$ is Ulrich we also get $N$ is Ulrich. If there is an irreducible morphism from $M \rt N$ then by our assumptions
on the ring there is also an irreducible morphism from $N \rt M$. By our earlier argument we get $N$ is Ulrich.

As there is no irreducible map from $A \rt M$ or from $M \rt A$ it follows that $\mathcal{U}$ defines a union of connected
components of $\Gamma(A)$.  
\end{proof}

\begin{remark}\label{odd-dim-ulrich}
 If $A = Q/(f)$ is a hypersurace ring with $\dim A$ is odd then note that Auslander-Reiten translate $\tau(M) = \Syz_1^A(M)$ which
 in general is not an Ulrich module (even if $M$ is Ulrich), see \cite[Theorem 2]{P}. So for odd dimensions our technique to produce
 an infinite family of indecomosable Ulrich modules with unbounded multiplicities; fails.
\end{remark}

\section{Proof of Theorem \ref{rev}}

In this section we give a proof of Theorem \ref{rev}.

\s We recall the definition of linkage of modules as given in \cite{MS}.  Throughout 
$(A,\m)$ is  a Gorenstein local ring of dimension $d$.

\s Let us recall the definition of transpose of a module. Let $F_1 \xrightarrow{\phi} F_0 \rt M \rt 0$ be a minimal presentation of $M$.  Let $(-)^* = \Hom(-,A)$. The \textit{transpose} $\Tr(M)$ is defined by the exact sequence
\[
0 \rt M^* \rt F_0^* \xrightarrow{\phi^*} F_1^* \rt \Tr(M) \rt 0.
\] 
Set $\lambda(M) = \Syz_1(\Tr(M))$.  We note that if $M$ is a maximal \CM \ $A$-module then $\Tr(M) = (\Syz_2(M))^*$.
\begin{definition}
Two $A$-modules $M$ and $N$ are said to be \textit{horizontally linked} if 
$M \cong \lambda(N)$ and $N \cong \lambda(M)$.
\end{definition}
If $E$ is a stable maximal \CM \ $A$-module then it is known  that $E$ is linked to $\lambda(E)$, i.e., $\lambda^2(E) = E$ see \cite[Corollary 7] {MS}.
 Note if  $M$ is an indecomposable  non-free maximal \CM \ $A$-module then so is $\lambda(M)$.
\begin{proof}
We prove the result only for $\lambda$. The proof for $D$ is in fact simpler. 

Let $M, N$ be indecomposable  non-free maximal \CM \ $A$-modules.  Using terminology from \ref{irr} it suffices to  prove that there exists an isomorphism
\[
\phi \colon \frac{(M,N)_1}{(M,N)_2} \rt \frac{(\lambda(N), \lambda(M))_1}{(\lambda(N), \lambda(M))_2}.
\] 
Let $f \in (M,N)_1$. Then as $M,N$ are indecomposable we get that  $f$ is not an isomorphism. In particular it is not a split mono. Let $f_2 \colon \Syz_2(M) \rt \Syz_2(M)$ be a lift of $f$. By \ref{lift-mono} and \ref{dual-irr} it follows that 
$f_2$ is not a split mono. Let $f_2^* \colon \Tr(N) \rt \Tr(M)$ be the dual of $f_2$. Then $f_2^*$ is not a split epi. Let $g_f \colon \lambda(N) \rt \lambda(M)$ be any lift of $f_2^*$. 
Define
\begin{align*}
\widetilde{\phi} \colon (M, N)_1  &\rt \frac{(\lambda(N), \lambda(M))_1}{(\lambda(N), \lambda(M))_2}, \\
f &\mapsto g_f + (\lambda(N), \lambda(M))_2.
\end{align*}
We first show that this map is independent of the choices we made. If $f_2^\prime$ is another lift of $f$ then $f_2 - f_2^\prime \in \beta(\Syz_2(M),\Syz_2(N))$. So
$f_2^* - (f_2^\prime)^* \in \beta(\Tr(N), \Tr(M))$. We know that if $\sigma \in \beta(\Tr(N), \Tr(M))$ then any lift of $\sigma$ is in 
$\beta(\lambda(N), \lambda(M))$.
Thus we have $g_f - g_f^\prime  = \delta \in \beta(\lambda(N), \lambda(M))$.
As $\lambda(N), \lambda(M)$ are indecomposable and non-free we get that 
$\delta \in (\lambda(N), \lambda(M))_2$. Thus $\widetilde{\phi}$ is well-defined.
It is elementary to show that $\widetilde{\phi}$ is $A$-linear.

Now let $f \in (M, N)_2$. Then there exists a maximal \CM \ $A$-module $X$ and 
 a commutative diagram
\[
\xymatrix{ 
M
\ar@{->}[d]_{u} 
\ar@{->}[dr]^{f}
 \\ 
X 
\ar@{->}[r]_{v} 
&N
}
\]
such that $u$ is not a split mono and $v$ is not a split epi.  Let
$u_2 \colon \Syz_2(M) \rt \Syz_2(X)$ be a lift of $u$ and 
 $v_2 \colon \Syz_2(X) \rt \Syz_2(N)$ be a lift of $v$.
By \ref{lift-mono} and \ref{lift-epi} we get that $u_2$ is not a split mono and $v_2$ is not a split epi. 
  Then $f_2 = v_2 \circ u_2$ is a lift of $f$. Then $f_2^* = u_2^* \circ v_2^*$.
  Also $u_2^*$ is not a split epi and $v_2^*$ is not a split mono. Let $\Syz_1(u_2^*)$ be a lift of $u_2^*$ and $\Syz_1(v_2^*)$ be a lift of $v_2^*$. Then
  $g_f = \Syz_1(u_2^*) \circ \Syz_1(v_2^*)$ is a lift of $f_2^*$.
  By \ref{lift-mono} and \ref{lift-epi} we get that $\Syz_1(u_2^*)$ is not a split epi and $ \Syz_1(v_2^*) $ is not a split mono. So $g_f \in (\lambda(N), \lambda(M))_2$. Thus we have a well-defined $A$-linear map
  \[
\phi \colon \frac{(M, N)_1}{(M, N)_2}  \rt \frac{(\lambda(N), \lambda(M))_1}{(\lambda(N), \lambda(M))_2}. 
\]
  As $\lambda^2(M) = M$ and $\lambda^2(N) = N$ we have a well defined  $A$-linear map
  \[
  \psi \colon  \frac{(\lambda(N), \lambda(M))_1}{(\lambda(N), \lambda(M))_2}  \rt \frac{(M, N)_1}{(M, N)_2}. 
  \]
  Finally it is tautological that $\phi$ and $\psi$ are inverses of each other. 
  Thus $\lambda \colon \Gau(A) \rt \Gau(A)^{rev}$ is an isomorphism.
  
Now assume that $A$ is not a hypersurface ring. We first note that
$\Syz_1(\lambda(M)) = M^*$ when $M$ is stable maximal \CM \ $A$-module. If $\lambda(M) = D(M)$ for all indecomposable maximal non-free $A$-modules then $M^*$ (and so $M$) has a periodic resilution with period $1$. It follows that $A$ is a hypersurface ring, a contradiction.

Next we show that there exist's $E$ with $D(E) \neq E$.  As $A$ is not a hypersurface there exists an MCM module $M$ which is not periodic. Let $M_1 = \Syz_1(M)$. As $M$ is not periodic either $M \neq M^*$ or $M_1 \neq M_1^*$.
  
  If $\lambda(M) = M$ for all indecomposable maximal \CM \ non-free $M$ then note that $\Syz_1(M) = M^*$ for all such $M$. We now note that
  \[
  \Syz_{-2}(M^*) \cong (\Syz_2(M))^* \cong \Syz_1(\Syz_2(M)) = \Syz_3(M)
  \]
  We now note that
  \[
 \Syz_{-2}(M^*) = \Syz_{-2}(\Syz_1(M)) =   \Syz_{-1}(M)
  \]
  It follows that $M$ is periodic for all indecomposable maximal \CM \ non-free $M$. Thus $A$ is a hypersurface, a contradiction.
  
  Thus $\lambda, D \colon \Gau(A) \rt \Gau(A)^{rev}$ are distinct isomorphism's if $A$ is not a hypersurface ring.  Furthermore $\lambda \neq 1$ and $D \neq 1$.
\end{proof}

The following result is immediate:
\begin{corollary}\label{syz}
(with hypotheses as in \ref{hypothesis}). For all $n \in \ZZ$ the map
$\Syz_n \colon \Gau(A) \rt \Gau(A)$ is an isomorphism.
\end{corollary}
\begin{proof}
We have $\Syz_1 \circ \lambda = D$. So $\Syz_1 = D \circ \lambda$ and $\Syz_{-1} = \lambda \circ D$. The result follows.
\end{proof}
We now give
\begin{proof}[Proof of Theorem \ref{syzn}.]  Let $[M]$ in $\mathcal{C}$. Recall
\[
 I(M) = \{ n \mid [\Syz_n(M)] \in \mathcal{C} \}.
\]
We first show that $I(M)$ is an ideal in $\ZZ$. As $M = \Syz_0(M)$ we get $0 \in I(M)$.
Now let $n \in I(M)$. The isomorphism $\Syz_{-n} \colon \Gau(A) \rt \Gau(A)$ maps $\mathcal{C}$ to itself since
$\Syz_{-n}(\Syz_n(M)) = M$. In particular we have $[\Syz_{-n}(M)] = \Syz_{-n}([M]) \in \mathcal{C}$.
 If $m,n \in \mathcal{C}$ then note that the isomorphism $\Syz_n \colon \Gau(A) \rt \Gau(A)$ maps $\mathcal{C}$ to itself
 as $\Syz_n([M]) = [\Syz_n(M)] \in \mathcal{C}$. Therefore $[\Syz_{n+m}(M)] = \Syz_n([\Syz_m(M)]) \in \mathcal{C}$. Thus
 $I(M)$ is an ideal in $\mathbb{Z}$. In particular there exist's a unique non-negative integer $i(M)$ such that
 $I(M) = i(M)\ZZ$.
 
 To prove rest  of the assertion of the theorem we first make a convention:  if $[X], [Y] \in \mathbb{C}$ then write
 $[X] <--> [Y]$ if there is an irreducible map from $X$ to $Y$ OR there is an irreducible map from $[Y]$ to $[X]$.
 
 As $[M], [N]$ are in $\mathcal{C}$ there is a sequence
 \[
  [M = X_0] <--> [X_1]  <--> \cdots <--> [X_{n-1}] <--> [X_n = N],
 \]
in $\mathcal{C}$. Set $a = i(M)$ and $b = i(N)$. By \ref{syz} we have the following sequence in $\Gau(A)$:
\[
 [\Syz_a(M)] <--> [\Syz_a(X_1)] <--> \cdots <--> [\Syz_a(N)].
\]
As $[\Syz_a(M)] \in \mathcal{C}$ we get $[\Syz_a(N)] \in \mathcal{C}$. So $a \in I(N)$ and therefore $I(M) \subseteq I(N)$.
Similarly we get $I(N) \subseteq I(M)$. Thus $I(M) = I(N)$.
\end{proof}

We now give

\begin{proof}[Proof of Corollary \ref{conn-0}]
(1) Suppose if possible $\mathcal{D}$ has only finitely many vertices. Then $\Syz_n(\Dc)$ cannot be a component of $\widetilde{\Gamma}(A) = \Gamma(A) \setminus \Gamma_0(A)$. As $\Gau_0(A)$ has only finitely many components we get 
$\Syz_n(\Dc) = \Syz_m(\Dc)$ for some $n > m$. Set $c = n - m$. Then $\Syz_c(\Dc) = \Dc$. Therefore $\Syz_{lc}(\Dc) = \Dc$ for all $l \in \ZZ$.

We note that the function $\Syz_{lc}$ permutes vertices of $\Dc$ among itself. As $\Dc$ is finite it follows that all modules in $\Dc$ is periodic.

As $\Dc$ is a connected component of $\Gau_0(A)$ it follows that there exists
$[M] \in \Dc$ such that there is an irreducible map either from $M$ to $A$ or an irreducible map from $A$ to $M$. In 
the first case $M$ is a component of $X(\m)$ the maximal \CM \ approximation of $\m$. By \ref{mcm-extremal} we get that $M$ is extremal. As $M$ is periodic we get that $A$ is a hypersurface ring, a contradiction.
In the second case there is an irreducible map from $\Syz_{-d + 2}(M) \rt A$. Note 
as $M$ is periodic then so is 
$\Syz_{-d + 2}(M)$. An argument similar to the earlier case yields that $A$ is a hypersurface ring, a contradiction.

(2) Suppose if possible the function $f \colon Vert(\Dc) \rt \ZZ$ given by $f([M]) = e(M)$ is bounded.  
As $e(M) \geq \mu(M)$ and $e(\Syz_1(M)) = e(A)\mu(M) - e(M)$, it follows that the multiplicity function on $Vert(\Syz_1(\Dc))$ is 
bounded. Iterating we get that the multiplcity function on $Vert(\Syz_n(\Dc))$ is bounded for each $n \geq 1$.
Then $\Syz_n(\Dc)$ cannot be a component of $\widetilde{\Gamma}(A) = \Gamma(A) \setminus \Gamma_0(A)$. As $\Gau_0(A)$ has only finitely many components we get 
$\Syz_n(\Dc) = \Syz_m(\Dc)$ for some $n > m$. Set $c = n - m$. Then $\Syz_c(\Dc) = \Dc$. Therefore $\Syz_{lc}(\Dc) = \Dc$ for all $l \in \ZZ$. In particular there exists $c$ such that 
$$(**) \quad \quad  \beta_{il}(M) \leq c \ \quad \text{ for all } \ l \geq 0 \ \text{ and all $[M]$ in $\Dc$}.$$

As $\Dc$ is a connected component of $\Gau_0(A)$ it follows that there exists $[M] \in \Dc$ such that there
is an irreducible map either from $M$ to $A$ or an irreducible map from $A$ to $M$. 
In the first case $M$ is a component of $X(\m)$ the maximal \CM \ approximation of $\m$.
By \ref{mcm-extremal} we get that $M$ is extremal. As $A$ is not an hypersurface we get the following
\begin{enumerate}
\item
If $A$ is a complete intersection of codimension $c \geq 2$ then $\cx(M) = c$. Furthermore $\lim \beta_i(M) = \infty$.  In particular the sequence $\{\beta_{il}(M) \}$ is unbounded. Thus (**) is not possible in this case.
\item
If $A$ is Gorenstein but not a complete intersection then  \\ $\curv(M) = \curv(k) > 1$. So there exists $r > 1$ such that $\beta_i(M) > r^i$ for all $i \gg 0$. 
Thus (**) is not possible in this case too.
\end{enumerate}

In the second case note that there is an irreducible map from $N = \Syz_{ -d + 2}(M) $ to $A$. We then have that for all $i \geq 0$ 
\[
 \beta_{il + d - 2}(N) \leq c 
\]
Then an argument similar to above gives a contradiction. 
\end{proof}

\section{Obstruction to quasi AR-sequences}
Let the setup be as in \ref{hypothesis}. Let $M$ be a non-free maximal \CM \ indecomosable $A$-module.

\s \label{non-exist-qar} Let $s \colon 0 \rt \tau(M) \rt E_M \rt M \rt 0$ be the AR-sequence ending at $M$. Then using  Proposition \ref{ar-qar} we get 
the following:\\
There is no quasi-AR sequence ending at $M \iff$ $E_M$ is free.

The next result gives an essential obstruction to non-existence of quasi AR-sequences.
\begin{lemma}\label{no-q}
 [with hypothesis as in \ref{hypothesis}] Further assume $d  \neq 1$. Suppose there is a non-free indecomposable maximal
 \CM \ module $M$ such that there  is no quasi AR-sequence ending at $M$. Then $A$ is a hypersurface ring. 
\end{lemma}
\begin{proof}
 By \ref{non-exist-qar} it follows that $\tau(M) = \Syz_1(M)$. By construction \\  $\tau(M) = \Syz_{-d +2}(M)$. Therefore we get that
 $M \cong \Syz_{-d + 1}(M)$. As $d \neq 1$ we get that $M$ (and so $\tau(M)$) is periodic.
 
 $E_M$ is non-zero and free. In particular it has $A$ as a summand.
 So there is an irreducible map from $\tau(M) \rt A$. It follows that $\tau(M)$ is a summand of $X(\m)$, the maximal 
 \CM \ approximation of $\m$. By \ref{mcm-extremal} we get that $\tau(M)$ is extremal $A$-module. It is also periodic. 
 So $A$ is a hypersurface ring.
\end{proof}

We now analyze hypersuface rings having perhaps modules $M$ such that there is no quasi AR-sequence ending at $M$.
Let $<M>$ denote the isomorphism class of a module $M$.
Set
\[
 Q^c(A) = \{ <M> \mid [M] \in \Gau(A) \ \text{with no quasi AR-sequence ending at $M$}\}.
\]
We show
\begin{proposition}\label{c-no-q}
 Let $(A,\m)$ be a complete equicharacteristic hypersurface isolated singularity. Assume $d = \dim A$ is even and non-zero.
 Also assume that $k = A/\m$ is algebraically closed.
 Then
 \begin{enumerate}[\rm (1)]
  \item 
  $Q^c(A)$ is a finite set (possibly empty).
  \item
  If $A$ is not of finite representation type and $\Syz_d(k)$ is indecomposable then $Q^c(A)$ is empty.
  \item
  If $Q^c(A)$ is non-empty and if $<M> \  \in Q_c(A)$ then
  \begin{enumerate}[\rm (a)]
   \item 
   $\Syz_n(M) = M$ for all $n \in \ZZ$.
   \item
   $[M]$ is an isolated component of $\Gau(A)$.
  \end{enumerate}
 \end{enumerate}
\end{proposition}
\begin{proof}
 We note that as $A$ is a hypersurface and $d$ is even we get that $\tau(M) = M$ for any non-free maximal \CM \ indecomosable
 $A$-module $M$.
 
 (1) If $<M> \in Q^C(A)$ then there is an irreducibe map from $M = \tau(M) \rt A$. So $M$ is a component of $X(\m)$. It follows
 that $Q^c(A)$ is a finite set.
 
 (2) As $\Syz_d(k)$ is indecomposable there is a unique non-free component of $X(\m)$. It follows that $\sharp Q^c(A) \leq 1$.
 If $<M> \in Q^c(A)$ then note that $[M] \leftrightarrows [A]$ is a connected component of $\Gamma(A)$. It follows that
 $A$ is of finite representation type, a contradiction.

 (3)(a) As there is no quasi AR-sequence ending at $M$ we get that $E_M$ is free. So $\tau(M) = \Syz_1(M)$. As $\dim A$ is even we
 get $M = \tau(M)$. As $A$ is a hypersurface we get $\Syz_n(M) = M$ for all $n \in \ZZ$.
 
 (3)(b) Notice $[M]$ is only connected to $[A]$. So  we get that $[M]$ is an isolated component in $\Gau(A)$.
 \end{proof}

\section{Structure of $\Gamma_0(A)$}
In this section we completely determine the structure of $\Gamma_0(A)$ when $\dim A = 2$ and its multiplcity
$e(A) \geq 3$.
\begin{theorem}\label{det-conn-comp-2}
[with hypothesis as in \ref{hypothesis}.]  Assume $\dim A = 2$ and $e(A) \geq 3$. Set $M_1 = A$.
 Then $\Gamma_0(A)$ is of the form
\[
M_1 \leftrightarrows M_2 \leftrightarrows M_3 \leftrightarrows  M_4 \leftrightarrows \cdots  \leftrightarrows M_n  \leftrightarrows \cdots
\]
where $e(M_n) = ne(M_1)$ for all $n \geq 1$.  Furthermore 
\begin{enumerate}[\rm(1)]
\item
$X(\m) = M_2 \oplus F$ where $F$ is free.
\item
$M_n^* = M_n$ for all $n \geq 1$.
\end{enumerate}
\end{theorem}
\begin{remark}
 We do not have any idea of the structure of $\Gamma_0(A)$ when $e(A) = 2$ (so necessarily $A$ is a hypersurface)  and $A$ is of infinite
 representation type. The reason is that Proposition \ref{rank} given below breaks down in the case $e(A) = 2$.
\end{remark}

The following result is essential in our  proof of Theorem \ref{det-conn-comp-2}.
\begin{proposition}\label{rank}
(with hypotheses as in  Theorem \ref{det-conn-comp-2}). Let $X(\m)$ be a MCM approximation of $\m$. Write $X(\m) = M \oplus F$ where $F$ is free and $M$ has no free summands. Then
\begin{enumerate}[\rm (1)]
\item
$M$ is indecomposable.
\item
$\rank M  = 2$.
\item 
$M \cong M^*$.
\end{enumerate}
\end{proposition}
\begin{proof}
(1) By \cite[Theorem B]{T}; $\Syz_2^A(k)$ is indecomposable. So $X(k) = \Syz_2^A(k)^*$ is indecomposable.
As $X(\m) = \Syz_1^A(X(k)) \oplus G$ where $G$ is free we get that 
$M = \Syz_1^A(\Syz_2^A(k)^*)$ is indecomposable by \cite[8.17]{Y}.

(2) We get $M^* = \Syz_{-1}^A(\Syz_2^A(k))$. Let $x,y$ be a $A \oplus M \oplus \Syz_2^A(k)$-superficial sequence.
Set $C = A/(x,y)$. If $E$ is an $A$-module then set $\ov{E} = E/(x,y)E$.  Notice
\[
\ov{M^*} \cong \Syz_{-1}^{C}(\ov{\Syz_2^A(k)}) \quad \text{and} \quad \ov{\Syz_2^A(k)} \cong 
\Syz_2^C(k)\oplus \Syz_1^C(k)^2\oplus \Syz_0^C(k). 
\]
Therefore
\[
\ov{M^*}  \cong \Syz_1^C(k)\oplus \Syz_0^C(k)^2\oplus \Syz_{-1}^C(k). 
\]
We note that as we have an exact sequence $0 \rt k = \soc(C) \rt C \rt C/\soc(C) \rt 0$. Thus $\Syz_{-1}^C(k) = C/\soc(C)$. Let $\n$ be the maximal ideal of $C$. Thus we have
\[
\ov{M^*}  \cong  \n \oplus k^2 \oplus C/\soc(C).
\]
Thus $\ell(\ov{M^*}) = 2\ell(C)$. Therefore
\[
e(M) = e(M^*) = e(\ov{M^*}) = \ell(\ov{M^*}) = 2\ell(C) = 2 e(C) = 2e(A).
\]
It follows that $\rank M =  2$.

(3) As there is an irreducible map $M \rt A$ there exists an irreducible map $A \rt M^*$.
As $\dim A = 2$ we have $\tau(M^*) = M^*$. So there is an irreducible map from 
$M^* \rt A$. Thus $M^*$ is an irreducible component of $X(\m)$. By (1) we have $M^* \cong M$.
\end{proof}

We now give
\begin{proof}[Proof of Theorem \ref{det-conn-comp-2}.] Set $X(\m) = M_2 \oplus F$ where $F$ is free and $M_2$ has no free summands. By Proposition \ref{rank} we get that $M_2$ is indecomposable of rank 2. We have the AR-sequence
\[
0 \rt M_2 \rt M_1^a \oplus X \rt M_2 \rt 0.
\]
Thus $a + \rank X = 4$. By Lemma \ref{no-q} and Proposition \ref{c-no-q}(2) we
get that $X \neq 0$.  Thus $1 \leq a \leq 3$.
We assert $a = 1$. We prove this by showing that the cases $a = 2$ or $3$ do not occur. \\
Claim 1: $a \neq 3$.\\
Suppose if possible $a = 3$ then rank $X$ is one. So $X$ is indecomposable. As $\dim A = 2$ and there is an irreducible map from $X$ to $M_2$, there is an irreducible map from $M_2 \rt X$. By rank considerations we get that the AR-quiver ending at $X$ is
\[
0 \rt X \rt M_2  \rt X \rt 0.
\] 
It follows that $M_1, M_2$ and $X$ constitute a connected component of $\Gamma_0(A)$ and so it is equal to $\Gamma_0(A)$.
Therefore $A$ has finite representation type, a contradiction.

Claim 2: $a \neq 2$.\\
If possible assume $a = 2$.
It follows that $\rank X = 2$. We  assert:\\
Subclaim 3: $X$ is indecomposable.\\
Suppose if possible $X = X_1 \oplus X_2$ where $\rank X_i = 1$. 
As $\dim A = 2$ and there is an irreducible map from $X_i$ to $M_2$, there is an irreducible map from $M_2 \rt X_i$. By rank considerations we get that the AR-quiver ending at $X_i$  for $i = 1, 2$ is
\[
0 \rt X_i \rt M_2  \rt X _i\rt 0.
\] 
It follows that $M_1, M_2$, $X_1$ and  $X_2$ constitute a connected component of $\Gamma_0(A)$ and so it is equal to $\Gamma_0(A)$. It follows that $A$ has finite representation type, a contradiction. Thus $X$ is indecomposable.

The AR-sequence ending at $X$ is
\[
0 \rt X \rt M_2 \oplus X_1 \rt X \rt 0.
\]
By an argument similar to Subclaim-3 we get that $X_1$ is indecomposable of rank 2. Set $X_0 = X$. 

For $i \geq 1$, by   an argument similar to Subclaim-3 we get that there exists indecomposable module $X_{i+1}$ of rank 2 such that the AR-sequence ending at $X_i$ is
\[
0 \rt X_i \rt X_{i-1} \oplus X_{i+1} \rt X_i \rt 0.
\]
Thus $\Gamma_0(A)$ consists of the modules $\{ M_1, M_2, X_i \mid i \geq 0 \},$ Also rank $X_i = 2$. This implies that $A$ is  of finite representation type (see \cite[6.2]{Y}), a contradiction. 

By claims 1, 2 we get $a = 1$. Thus $\rank X = 3$. \\
Claim 4: $X$ is indecomposable. \\
Suppose if possible this is not so. Then either \\
Subcase 5: $X = X_1 \oplus X_2 \oplus X_3$ where $\rank X_i = 1$ for $1\leq i \leq 3$, OR \\
Subcase 6: $X = X_1 \oplus X_2 $ where $\rank X_i = i$ for $i = 1,2$.\\
We show that subcase 5, 6 are not possible. If subcase 5 occurs then by rank considerations the AR-quiver ending at $X_i$ is
\[
0 \rt X_i \rt M_2 \rt X_i \rt 0 \quad \text{for} \ i = 1,2,3.
\]
Thus the vertices of $\Gamma_0(A)$ will be 
\[
\{M_1, M_2, X_1, X_2, X_3 \}.
\]
This implies that $A$ has finite representation type, a contradiction.

If subcase 6 occurs 
then by rank considerations the AR-quiver ending at $X_1$ is
\[
0 \rt X_1 \rt M_2 \rt X_1 \rt 0. 
\]
Furthermore the AR-quiver ending at $X_2$ is
\[
0 \rt X_2 \rt M_2 \oplus X_3 \rt X_2 \rt 0. 
\]
Note $\rank X_3 = 2$. By an argument similar to that of subcase 5 we get that $X_3$ is indecomposable. Iterating  we obtain rank two indecomposable modules $X_i$ for $i \geq 4$ such that  the AR-quiver ending at $X_i$ is
\[
0 \rt X_i \rt X_{i+1} \oplus X_{i-1} \rt X_i \rt 0. 
\]
It follows that the vertices of $\Gamma_0(A)$ is 
\[
\{ M_1, M_2, X_i \mid i \geq 1 \}.
\]
As there is a bound on the ranks of vertices of $\Gamma_0(A)$ it follows that $A$ is of finite representation type, a contradiction.

 Set $M_3 = X$.  We have rank $M_3 = 3$ and that $M_3$ is indecomposable. Inductively assume that we have indecomposable MCM $A$-modules $M_1,\ldots, M_n$ with $n \geq 3$ and 
 $\rank M_i = i$ such that the AR-sequence ending at $M_j$ for $j \leq n-1$ is
 \[
 0 \rt M_j \rt M_j \oplus M_{j+1} \rt M_j \rt 0.
 \]
 Let the AR-sequence ending at $M_n$ be 
 \[
 0 \rt M_n \rt M_{n-1} \oplus Y \rt M_n \rt 0.
 \]
 Clearly $\rank Y = n+1$. If we prove that $Y$ is indecomposable then we can set $M_{n+1} = Y$ and we will be done by induction.
 
 Let $Z$ be an indecomposable summand of $Y$. Then the AR-sequence ending at $Z$ is
 \[
 0 \rt Z \rt M_n \oplus W \rt Z \rt 0,
 \]
 where $W$ is an MCM $A$-module (possibly zero). Nevertheless we get that $\rank Z \geq n/2$.
 
 As $n \geq 3$,  $\rank Y  = n +1$  and an indecomposable summand $Z$ of $Y$ has rank atleast $n/2$ it follows that $Y$ has at most two indecomposable summnads. 
 
 We want to prove that $Y$ is indecomposable. Suppose it is not so. Then by our previous argument it has two indecomposable summands say $Y_1$ and $Y_2$. Suppose $\rank Y_1 \leq \rank Y_2$. Then we have
 \[
 \frac{n}{2} \leq \rank Y_1 \leq \frac{n+1}{2}.
 \]
 
 We consider two cases:\\
 Case 1: $n = 2m + 1$ is odd.\\
 We get $\rank Y_1 = m + 1$.  So $\rank Y_2 = m + 1$ also. Let the AR-sequence ending at $Y_1$ be
 \[
 0 \rt Y_1 \rt M_n \oplus T \rt Y_1 \rt 0.
 \]
 Thus $T$ has rank $1$. The AR-sequence ending at $T$ is
 \[
 0 \rt T \rt Y_1 \oplus L \rt T \rt 0.
 \]
 As $m + 1 \leq 2$ we get $m \leq 1$. As $m \geq 1$ we get $m = 1$. Therefore $n = 2m+1 = 3$. Now consider the case $n = 3$. We get $\rank Y_j = 2$ for $j = 1, 2$ and $\rank T = 1$. Furthermore $L = 0$. Similarly the AR-sequence ending at $Y_2$ will be
 \[
 0 \rt Y_2 \rt M_3 \oplus T^\prime \rt Y_2 \rt 0,
 \]
 where $T^\prime$ has rank $1$. The AR-sequence ending at $T^\prime$ is
 \[
 0 \rt T^\prime \rt Y_2  \rt T^\prime \rt 0.
 \]
 It follows that the vertices of $\Gamma_0(A)$ will be 
 \[
 \{ M_1, M_2, M_3, Y_1, Y_2, T, T^\prime \}.
 \]
 It follows  that $A$ has finite representation type, a contradiction.
 
 Case 2: $n = 2m$ is even.\\
 We get $\rank Y_1 = m$ and $\rank Y_2 = m + 1$.
 The AR sequence ending at $Y_1$ is
 \[
 0 \rt Y_1 \rt M_n \rt Y_1 \rt 0.
 \]
 The AR sequence ending at $Y_2$ is
 \[
 0 \rt Y_2 \rt M_n \oplus T \rt Y_2 \rt 0.
 \]
 It follows that rank $T = 2$. We have to consider two sub cases: \\
 Subcase-1: $T$ is decomposable. In this case $T = T_1 \oplus T_2$ where $\rank T_i = 1$ for $i = 1, 2$. The AR-sequence ending at $T_1$ is
 \[
  0 \rt T_1 \rt Y_2 \oplus L \rt T_1 \rt 0.
 \]
 We have $2 = m+1 + \rank L$. As $m \geq 1$ we get $m = 1$ and $L = 0$. So $n = 2$. We have already dealt with this case.
 
 Subcase-2: $T$ is indecomposable. The AR-sequence ending at $T$ is
 \[
 0 \rt T \rt Y_2 \oplus W \rt T \rt 0.
 \] 
 We have $ 4 = m+1 + \rank W$. As $m \geq 1$ the possibilities for $m$ is $1, 2, 3$.
 If $m = 1$ then $n = 2$. This case has been discussed earlier. Next we consider the case $m = 3$. In this case $W = 0$. So the vertices of $\Gamma_0(A)$ will be
 \[
 \{ M_i, Y_1, Y_2, T \mid 1 \leq i \leq n \}.
 \]
It follows  that $A$ has finite representation type, a contradiction.
 
 Finally we consider the case when  $m = 2$. So $n = 4$.  Thus $ \rank W = 1$. 
 The AR-sequence ending at $W$ is 
 \[
 0 \rt W \rt T \rt W \rt 0.
 \]
 Thus the vertices of $\Gamma_0(A)$ will be
 \[
 \{ M_i, Y_1, Y_2, T, W \mid 1 \leq i \leq n \}.
 \]
 It follows that $A$ has finite representation type, a contradiction.
 
 (2) We note that the dual map $D \colon \Gau(A) \rt \Gau(A)^{rev}$ is an isomorphism of graphs. As $D(M_2) = M_2^* \cong M_2$ and as $\Gau_0(A)$ is connected we get that $D$ maps $\Gau_0(A)$ to itself. Comparing ranks we get $M_n^* \cong M_n$ for all $n \geq 3$.
\end{proof}

\section{Proof of Theorem \ref{det-conn-comp} and  Corollary \ref{d2-conn}}
 In this section of the results as stated in the title of this section. Throughout $(A,\m)$ is an equicharacteristic Gorenstein 
 isolated singularity of dimension two. We also assume that $A$ is complete and the residue field $k$ is algebraically closed. Furthermore
 we assume that $e(A) \geq 3$.

 We first give 
\begin{proof}[Proof of Corollary \ref{d2-conn}]
 It suffices to show that $\Syz_n(M) \notin \Gau_0(A)$ for all $M \in \Gau_0(A)$ and for all $n \neq 0$. Using the terminology
 of Theorem \ref{syzn} we need to show $I(M) = 0$ for all $M$ in  
 $\Gau_0(A)$. We also recall that $I(M) = I(N)$ for all 
 $M, N \in \Gau_0(A)$. We denote this common value by $c$.
  
  We want to show $c = 0$. If possible assume $c> 0$. 
  Set
  \[
   V  = \{ |i-j| \mid M_j = \Syz_nM_i \ \text{for some} \ n \neq 0 \} \quad \text{and} \ r = \min V.
  \]
Notice $c \neq 0$ if and only if $V \neq \emptyset$. 

We first consider the case when $r = 0$. Say $M_i = \Syz_n M_i$ for some $n \neq 0$. We may assume $n > 0$. Then $M_i$ is periodic.
As $A$ is not a hypersurface this is a contradiction by Theorem \ref{first} and  Theorem \ref{mcm-extremal}.

We now assume $r \geq 1$. Say $M_{i + r} = \Syz_n(M_i)$ for some $r > 0$ and for some $n \neq 0$. Note we are not assuming
$n > 0$. As we have an irreducible morphism from $M_{i+r-1} \rt M_{i+r}$ we have an irreducible map from 
\[
 \Syz_{-n}(M_{i+r -1}) \rt M_i.
\]
So we have $M_{i+1} = \Syz_{-n}(M_{i+r -1})$ or $M_{i -1} = \Syz_{-n}(M_{i+r -1})$. The first case cannot occur as $r = \min V$.
So $M_{i -1} = \Syz_{-n}(M_{i+r -1})$ and therefore $M_{i+r-1} = \Syz_n(M_{i-1})$. Iterating this procedure we get that 
$M_{2+r} = \Syz_n(M_2)$.  We have irreducible maps from $M_{2+r-1}$ and $M_{2+r + 1}$ to $M_{2+r} = \Syz_n(M_2)$. So we have an irreducible map from $\Syz_{-n}(M_{2+r-1})$ and $\Syz_{-n}(M_{2+r+1})$ to $M_2$. It follows that atleast one of  $\Syz_{-n}(M_{2+r-1})$ and $\Syz_{-n}(M_{2+r+1})$ is $A$. This is a contradiction. 
\end{proof}

Next we give
 \begin{proof}[Proof of Theorem \ref{det-conn-comp}]
 The  assertion on the structure of $\mathcal{C}$ follows from Theorem \ref{det-conn-comp-2} and \cite[4.16.2]{Ben}.
  
  (1) This follows from Theorem \ref{det-conn-comp-2}.
  
  (2)(a) Let $\mathcal{C}$ be a connected component of $\Gamma(A)$ such that $[M] \in Vert(\mathcal{C})$ is a periodic module. Then
  by Theorem \ref{first} all the modules $N$ in $Vert(\mathcal{C})$ is periodic. We note that $\Syz_n(\mathcal{C})$ consists of periodic
  modules and so $[A] \notin Vert(\Syz_n(\mathcal{C})$ for all $n \in \ZZ$ (see Theorem \ref{mcm-extremal}). Using Theorem \ref{lift-qar}
  and Corollary \ref{qar-ar} we get that if $[M] \in Vert(\mathcal{C})$ and if $0 \rt M \rt E_M \rt M \rt 0$ is an AR sequence
  ending at $M$ then for all $n \in \ZZ$  the AR-sequence ending at $\Syz_n(M)$ is of the form $0 \rt \Syz_n(M) \rt \Syz_n(E) \rt
  \Syz_n(M) \rt 0$.
  
  Now consider the structure of $\mathcal{C}$ as given in (1).
  Let period of $M_1$ be $c$. We first show that $I(M_1) = c \ZZ$ (notation as in Theorem \ref{syzn}). Note $c \in I(M_1)$. 
  If $I(M_1) \neq c \ZZ$ then there exists $a$ with $1 < a < c$ such that $[\Syz_a(M_1)] \in Vert(\mathcal{C})$. We note that
  $0 \rt \Syz_a(M_1) \rt \Syz_a(M_2) \rt \Syz_a(M_1) \rt 0$ is the AR sequence ending at $\Syz_a(M_1)$. As $M_1$ is the unique 
  vertex in $\mathcal{C}$ which is connected to only one other vertex we get that $\Syz_a(M_1) = M_1$. This contradicts the fact that
  period of $M_1$ is $c$.
  
  We show by induction on $n \geq 2$ that the  period of $M_n$ is $c$. We first consider the case $n = 2$. As period of $M_1$ is $c$
  we get that $0 \rt M_1 \rt \Syz_c(M_2) \rt M_1 \rt 0$ is also an AR-sequence ending at $M_1$. By uniqueness of AR sequences we get
  $M_2 \cong \Syz_c(M_2)$.  Suppose for some $a$ with $1 \leq a < c$ we have $\Syz_a(M_2) = M_2$ then note that $a \in I(M_2) = I(M_1) =c\ZZ$, a contradiction.
  Thus period of $M_2$ is $c$.
  
  Now assume that period of $M_1, \ldots, M_n$ is $c$. We prove that period of $M_{n+1}$ is also $c$. As the period of $M_{n-1}$ and 
  $M_n$ is $c$ we get that $ 0 \rt M_n \rt M_{n-1}\oplus \Syz_c(M_{n+1}) \rt M_n \rt 0$ is another AR-sequence ending at $M_n$. 
  By uniqueness of AR-sequences we get that $M_{n+1} \cong \Syz_c(M_{n+1})$. Suppose for some $a$ with $1 \leq a < c$ we have
  $M_{n+1} = \Syz_{a}(M_{n+1})$. Then $a \in I(M_{n+1}) = I(M_1) = c \ZZ$, a contradiction. Thus period of $M_{n+1}$ is $c$. 
  The result follows.
  
  (2)(b) By \ref{d2-conn} there exists \emph{atmost} one $m_0 \geq 1$ such that $\Syz_{m_0}(\mathcal{C}) = \Gau_0(A)$. Thus for 
  $n > m_0$ we have that $[A] \notin Vert(\Syz_n(\mathcal{C}))$. Set $M_0 = 0$. We have that for all $n > m_0$ the sequence
  $0 \rt \Syz_n(M_i) \rt \Syz_n(M_{i-1}) \oplus \Syz_{n}(M_{i+1}) \rt \Syz_n(M_i) \rt 0$ is the AR quiver ending at $M_i$ for all
  $i \geq 1$. By Lemma \ref{min-gen} we get that for all $n > m$ and for all $i \geq 1$
  \[
   2 \beta_n(M_i) = \beta_n(M_{i-1}) + \beta_n(M_{i+1}).
  \]
As $M_0 = 0$ an easy recurssion yields that $\beta_n(M_i) = i \beta_n(M_1)$. The result follows.
 \end{proof}
 \section{curvature and complexity}
If $(A,\m)$ is a complete intersection of codimension $c$ then it is known that for any non-zero module $M$ we have
$0 \leq \cx M \leq c$. Furthermore for any integer $i$ with $0 \leq i \leq c$ there exists an $A$-module $M$ with 
complexity $i$. If $A$ is not a complete intersection then $\cx k = \infty$. To deal with this situation the notion of curvature
was introduced. It can be shown that $1 < \curv k  < \infty$ and for any non-zero module with infinite projective dimension
we have $1 \leq \curv M \leq \curv k$. Furthermore if $\cx M < \infty $ then $\curv M = 1$.   We first prove
\begin{proposition}\label{cx-curv}
 Let $(A,\m)$ be an equicharactersitic complete Gorenstein isolated singularity with algebraically closed residue field $k$. Assume 
 $A$ is not a complete intersection. Then
\begin{enumerate}[\rm (1)]
 \item 
 For any $i \geq 1$ the modules $M$ with complexity $i$ form a union of  connected components of $\Gamma(A)$.
 \item
 For any $\alpha \in [1, \curv k)$  the modules $M$ with curvature $\alpha$ form a
 union of  connected components of $\Gamma(A)$.
\end{enumerate}
\end{proposition}

We first show
\begin{lemma}\label{curv-lemma}[with hypotheses as in Proposition \ref{cx-curv}]
Let $1 \leq \alpha < \curv k$. Let $\mathcal{V_\alpha}$ be the collection of all indecomposable modules $M$ with $\curv M \leq \alpha$. Then $V_\alpha$ is a union of 
connected components of $\Gamma(A)$. Furthermore $\Gamma_0(A) \nsubseteq V_\alpha$.
\end{lemma}
\begin{proof}
Let $M \in V_\alpha$. Note that $\tau(M) = \Syz_{-d+2}(M) \in V_\alpha$. As $\alpha < \curv k$ it follows that there is no irreducible map from $M$ to $A$ or from $A$ to $M$, see \ref{mcm-extremal}.

Clearly $\Syz_n(M) \in V_{\alpha}$ for all $n \in \ZZ$. By a similar argument as before  there is no irreducible map from $\Syz_n(M)$ to $A$ or from $A$ to $\Syz_n(M)$ for all $n \in \ZZ$.

Let $0 \rt \tau(M) \rt E_M \rt M \rt 0$ be the AR-sequence ending at $M$. By 7.13 and 7.14 we get that
\begin{enumerate}
\item
$0 \rt \Syz_n(\tau(M)) \rt \Syz_n(E_M) \rt \Syz_n(M) \rt 0$ is the AR-sequence ending at $\Syz_n(M)$ for all $n \geq 0$.
\item
$\beta_n(E_M) = \beta_n(M) + \beta_n(\tau(M))$ for all $n \geq 0$.
\end{enumerate}
Thus we have $\curv(E) \leq \alpha$. If there is an irreducible map from $N$ to $M$ then $N$ is a factor of $E_M$ and so $\curv(N) \leq \curv(E) \leq c$. Thus $N \in V_\alpha$. In a similar fashion if there is an irreducible map from $M$ to $N$ then also $N \in V_\alpha$. Thus $V_\alpha$ is a union of connected components of $\Gamma(A)$.  Also clearly $\Gamma_0(A) \nsubseteq V_\alpha$.
\end{proof}
As an immediate consequence we get 
\begin{corollary}\label{less-a}
[with hypotheses as in Proposition \ref{cx-curv}]
Let $1 < \beta < \curv k$. Let $\mathcal{U_\beta}$ be the collection of all indecomposable modules $M$ with $\curv M < \beta$. Then $U_\beta$ is a union of 
connected components of $\Gamma(A)$. Furthermore $\Gamma_0(A) \nsubseteq  U_\beta$.
\end{corollary}
\begin{proof}
Let $1 = \alpha_1 < \alpha_2 < \cdots < \alpha_n < \alpha_{n+1} < \cdots$ be any strictly monotonically increasing sequence converging to $\beta$.
Notice
\[
U_\beta = \bigcup_{n \geq 1}V_{\alpha_n}
\]
The result now follows from Lemma  \ref{curv-lemma}.
\end{proof}
We now give 
\begin{proof}[Proof of Proposition \ref{cx-curv}]
We first prove (2). Let $C_\alpha = $ the collection of modules with complexity $\alpha$. Notice (with notation as in 
Lemma \ref{curv-lemma} and Corollary \ref{less-a}
\begin{enumerate}[\rm (a)]
\item
$C_1 = V_1$.
\item
For $1 < \alpha < \curv(k)$ we have $C_{\alpha} = V_\alpha - U_\alpha$.
\end{enumerate}
Thus (2) follows.

(1) This is similar to (2). We have to prove  results analogous  to Lemma \ref{curv-lemma} and Corollary \ref{less-a} first.
\end{proof}

We now give 
\begin{proof}[Proof of Theorem \ref{BNP}]
Suppose $A$  has a module $M$ with bounded betti-numbers but  not periodic. Then note that $A$ is not a complete intersection.
We note that a MCM $A$-module $M$ will have bounded betti-numbers if and only if $\cx(M) \leq 1$. By Proposition \ref{cx-curv}, $\Dc$ the collection of all such modules defines a union of connected components of $\Gamma(A)$. We note that modules $M$ having a periodic resolution will form a subset $\mathcal{C}$ of $\Dc$. By  Theorem \ref{first} we get that $\mathcal{C}$ is a union of connected components of $\Gamma(A)$. It follows that $\Dc \setminus \mathcal{C}$ is a union of connected components of $\Gamma(A)$. If $M$ is not periodic but has a bounded resolution then
$[M] \in \Dc \setminus \mathcal{C}$. The result follows.
\end{proof}

\end{document}